\documentclass[eqno,fleqn]{article}
\usepackage{amsmath}
\usepackage{amssymb}
\usepackage[shortlabels]{enumitem}
\usepackage[onehalfspacing]{setspace}
\usepackage[T1]{fontenc}
\usepackage{graphicx}
\usepackage{xcolor}

\newcommand{\bbbh}{\mathbb{H}} 
\newcommand{\bbbn}{\mathbb{N}} 
 
\newtheorem {overview}{Overview} 
\newtheorem{satz}{Satz}[section] 
\newtheorem{definition}[satz]{Definition}

\newtheorem{theorem}[satz]{Theorem}

\newtheorem{example}[satz]{Example}

\newtheorem{agree}[satz]{Agreement}
\newtheorem{consequ}[satz]{Consequences}
\counterwithout{equation}{section}

\begin{document}
\pagestyle{headings} 
\thispagestyle{empty}

\title{Abstract computation over first-order structures. Part I: Deterministic and non-deterministic BSS RAMs}
\author{Christine Ga\ss ner}

\begin{center}

{\Large Abstract Computation over First-Order Structures \vspace{0.2cm}}\\{\Large Part I: Deterministic and Non-Deterministic BSS RAMs}
\vspace{0.7cm}\\

{\bf Christine Ga\ss ner} \footnote{In Part III and Part IV, we present a generalization of the results discussed at the CCC in Kochel in 2015. My thanks go to Dieter Spreen, Ulrich Berger, and the other organizers of the conference in Kochel. Further results were discussed at the CCA 2015. I would like to thank the organizers of this conference and in particular Martin Ziegler and Akitoshi Kawamura.} {\vspace{0.1cm}\\University of Greifswald, Germany, 2025\\ gassnerc@uni-greifswald.de}\\\end{center}

\vspace{0.7cm}

\begin{abstract} 
Most ideas about what an algorithm is are very similar. Basic operations are used for transforming objects. The evaluation of internal and external states by relations has impact on the further process. A more precise definition can lead to a model of abstract computation over an arbitrary first-order structure. Formally, the algorithms can be determined by strings. Their meaning can be described purely mathematically by functions and relations derived from the operations and relations of a first-order structure. Our model includes models of computability and derivation systems from different areas of mathematics, logic, and computer science. To define the algorithms, we use so-called programs. Since we do this independently of their executability by computers, the so-called execution of our programs can be viewed as a form of abstract computation. This concept helps to highlight common features of algorithms that are independent of the underlying structures. Here, in Part I, we define BSS RAMs step by step. In Part II, we study Moschovakis' operator which is known from a general recursion theory over first-order structures. Later, we study hierarchies defined analogously to the arithmetical hierarchy by means of quantified formulas of an infinitary logic in this framework.

\vspace{0.6cm}


\end{abstract}

\vspace{0.5cm}

\newpage

\section{Introduction} 

Because many algorithms are used in different areas regardless of their executability by computers, we investigate algorithms in a very general framework that combines logic and theoretical computer science in a natural way. We present a concept of BSS RAMs whose development was influenced by several models of computation introduced in 1989 and later. Our machines over first-order structures help to describe sequences of computation steps and other derivation steps in a standardized way. We generalize various concepts of computation including the processing of strings or the codes of natural numbers by (deterministic and non-deterministic) Turing machines (cf.\,\cite{TURING37}) and the computation over the reals (cf.\,\cite{BSS89}) by using the BSS model for describing practical relevant algorithms on a high level of abstraction. This means that we consider a model of computation where the individual objects of the underlying universe are the smallest units that can be processed by means of basic operations. And, due to an input and an output procedure, it can be easier to handle all instances of one problem uniformly by a single BSS RAM. For this purpose, several features of the BSS machines\footnote{Note, that an important feature is missing from the description on the corresponding English Wikipedia website, January 25, 2025.} introduced by Lenore Blum, Michael Shub, and Stephen Smale in \cite{BSS89} were directly generalized and adjusted to better handle first-order structures that contain neither an infinite recursively enumerable subset nor a suitable set of constants nor any constant (as in \cite{GASS20}). This makes it possible, on the one hand, to describe the time complexity of a problem as a function whose values depend on which basic operations are permitted and, on the other hand, to obtain a common generalization of the complexity theories dealing with the computation over different first-order structures.

Here, in Part I, we gradually expand the possibilities to describe algorithms to show that the BSS RAMs and the use of oracles are a good tool for dealing with theoretical questions of abstract computability and the consequences. Here and later, we will deal with the following questions. Can we describe the essential properties of an algorithm, regardless of whether it can be executed by a computer or not? What properties do abstract machines need to provide a good characterization of our algorithmic reasoning at a level where the specific properties of the processed objects are largely disregarded? What are the consequences of allowing non-deterministic computational models? Which properties of the arithmetical hierarchy are retained when the use of an enumerable infinite set of individuals is not possible? We try to answer these questions.

Generally speaking, an algorithm is a process for solving or handling problems (cf.\,\,\cite{KleinEnz}). Starting from certain objects, which themselves or whose codes can be used as input, results can be returned after processing these objects. Various models --- such as the Turing machine \cite{TURING37}, non-machine-oriented models \cite{KLEENE,MOSCHO} and RAM models for computation over the natural numbers \cite{AHU74} or over the real numbers \cite{PS85} --- help to systematically describe the processing of numbers or other objects. More information about the great variety of models can be found in \cite{ ShS63, Hennie, S67, B92, HU79, WEIHRAUCH, BCSS98, H98, BH98}, and \cite{BP18}, to name just a few publications. The set of natural numbers or completely different first-order structures with their operations and relations can serve as the basis for processing objects. In connection with the introduction of structured programming, various machine-oriented models of computation were discussed in detail. The requirements for such models depend on the expectation that certain operations, tests, and transport commands should be executable where transporting an object means copying it by using auxiliary operations for manipulating memory addresses. The different models allow the discussion of complexity issues at a higher level of abstraction and thus help to obtain an initial estimate of the degree of difficulty of certain problems. On the other hand, this of course also raises many new theoretical questions, which have also promoted the design of the uniform BSS model \cite{BSS89} for the treatment of real numbers. 

We, in turn, have extended this model and introduced few further commands for the transport of arbitrary objects. This extension has the following consequences. It is not necessary to calculate addresses for copying objects by using the individuals of the underlying structures as codes of these addresses. This option is interesting, for example, if the structures, underlying the algorithmic processing, are only suitable for modeling biological processes as considered in \cite{Zerjatke} and not for ordinary calculations. Generalizations of the machine-oriented concepts such as presented in \cite{P95} and \cite{H98} were developed in the early 1990s. In \cite{GASS96} and \cite{GASS97}, the BSS model was directly generalized. The latter is in a certain sense close to a minimal variant (like the model in \cite{P95}), which nevertheless makes it possible to take a wide variety of structures into account. The number of types of instructions is small but sufficient. The further design was significantly influenced by the discussions in \cite{KLEENE,KOZEN,SOARE} and many others (see also\cite{BDGI,BlumN,Hoff09,Papad,Reischuk99,Schoening08,TuschikWolter}).

Based on the systematic introduction of the deterministic BSS RAMs in \cite{GASS20}, we extend this concept. 
 We will also consider so-called non-deterministic algorithms and describe such algorithms by means of several types of non-deterministic BSS RAMs. This should mean that we want to characterize algorithms that can be influenced in a non-deterministic way by options that make it possible to guess which individuals can be used additionally or which instructions could be selected to execute the next step, and the like. In Section \ref{DefBSSRAM}, we define different types of BSS RAMs step by step. Some examples are presented and discussed in Section \ref{SectionExamples}. Finally, we give a summary and an outlook.

In Parts II and III, we will discuss the simulation of Moschovakis' operator by non-deterministic BSS RAMs and the simulation of non-deterministic BSS RAMs by universal BSS RAMs in order to be able to justify the completeness of certain halting problems without the need to have constants. In Part IV, we will syntactically define a hierarchy of abstract decision problems and investigate the completeness of some problems for structures with a finite number of basic operations and relations, a semi-decidable identity relation, and a semi-decidable set of constants. In this way, it will be possible to characterize the syntactically defined hierarchy semantically by non-deterministic oracle machines. 

\section{BSS-RAMs over first-order structures}\label{DefBSSRAM}

Every considered algorithm can be syntactically described by means of a program that consists of a finite sequence of instructions some of which are given in Overview \ref{Sigma_InstructionsFini}. The instructions whose form depends on the used signature $\sigma$ can be interpreted by using a first-order structure ${\cal A}$ of this signature. Briefly and roughly speaking, an ${\cal A}$-machine that is able to add two numbers in one time unit should be able to execute instructions such as $A:= B+D$ in one single step. Here, $A:= B+D$ stands for the commands {\sf compute the sum $B+D$ and assign the sum to $A$}. The program of a corresponding BSS RAM could be 

\vspace{0.1cm}

\quad {\sf $1: \,$$ Z_1:= f_3^2(Z_2,Z_4);$\, $2: \,$ stop}.

\vspace{0.1cm}

\noindent where $Z_1$, $Z_2$, and $Z_4$ are {\em names} for so-called {\em $Z$-registers} used for storing numbers. The term $f_3^2(Z_2,Z_4)$ after the symbol $:=$ is given in prefix notation. $f_3^2$ is a symbol for a binary operation which we generally denote by $f_3$. 

\subsection{First-order structures and their signatures} 

Let us use any first-order structure ${\cal A}$ of the form $(U; (c_i)_{i\in N_1}; (f_i)_{i\in N_2}; (r_i)_{i\in N_3})$ for interpreting the descriptions of algorithms that are given in the form of suitable programs. Let $U$ be a non-empty set of individuals. For suitable sets $N_1$, $N_2$, and $N_3$, let $( c_i)_{i\in N_1}$ be a family of constants $c_i\in U$, $( f_i)_{i\in N_2}$ be a family of operations $ f_i:U^{m_i}\to U$ of arity $m_i$, and $( r_i)_{i\in N_3}$ be a family of relations $r_i$ of arity $k_i$ such that $r_i\subseteq U^{k_i}$. Let $U_{\cal A}$ be the universe $U$ of ${\cal A}$. We define ${\cal A}$-machines so that any ${\cal A}$-machine ${\cal M}$ will be able to execute its own program denoted by ${\cal P}_{\cal M}$. Since any program that we permit is a finite string, we need only a reduct ${\cal B}$ of the structure ${\cal A}$ in order to interpret and execute the program ${\cal P}_{\cal M}$. Such a reduct ${\cal B}$ can be a structure of the form $(U_{\cal A};c_{\alpha_1},\ldots,c_{\alpha_{n_1}};f_{\beta_1},\ldots, f_{\beta_{n_2}}; r_{\gamma_1},\ldots,r_{\gamma_{n_3}})$ if $\{\alpha_1,\ldots,\alpha_{n_1}\}\subseteq N_1$,$\{\beta_1,\ldots,\beta_{n_2}\}\subseteq N_2$, and $\{\gamma_1,\ldots,\gamma_{n_3}\}\subseteq N_3$ are satisfied \footnote{Sequences such as $c_{\alpha_1},\ldots,c_{\alpha_{0}}$ and $\alpha_1,\ldots,\alpha_{0}$ are empty and thus sequences of length $0$.} for some integers $n_1, n_2$, and $n_3$ in $ \bbbn=_{\rm df}\{0,1,2,\ldots\}$ with $n_1\leq |N_1|$, $n_2\leq |N_2|$, and $n_3\leq |N_3|$. The signature $(n_1;m_{\beta_1}, \dots,m_{\beta_{n_2}}; k_{\gamma_1}, \dots,k_{\gamma_{n_3}})$ of such a reduct is finite and can also be the signature of ${\cal P}_{\cal M}$. In Part III, we generally consider structures of the form $(U; (c_i)_{i\in N_1};f_1,\ldots, f_{n_2}; r_1,\ldots,r_{n_3})$. 

Thus, for simplicity, 
let now ${\cal A}$ be any first-order structure and $\sigma $ be any signature $(n_1;m_1, \dots,m_{n_2}; k_1, \dots,k_{n_3})$ with $n_1,n_2,n_3\in \bbbn$. We assume that the reduct ${\cal B}$ of ${\cal A}$ that a machine ${\cal M}$ over ${\cal A}$ can use to execute its program ${\cal P}_{\cal M}$ has the form $(U_{\cal A};c_{\alpha_1},\ldots,c_{\alpha_{n_1}};f_1,\ldots, f_{n_2}; r_1,\ldots,r_{n_3})$ and is of signature $\sigma$. 

\subsection{The $\sigma$-programs} 

Let ${\sf P}_{\sigma}$ be the set of all {\em$\sigma$-programs ${\cal P}$} that are strings of the form 
\vspace{0.1cm}

\quad 1: {\sf instruction}$_1$;\, 2: {\sf instruction}$_2$; \ldots;\, $\ell _{\cal P}-1$: {\sf instruction}$_{\ell _{\cal P}-1}$;\, $\ell _{\cal P}$: {\sf stop}.

\vspace{0.1cm}

\noindent Each substring ${\sf instruction}_\ell$ in ${\cal P}$ is a so-called {\em instruction} whose {\em label} is $\ell$ and whose form depends in part on the signature $\sigma$. We will here distinguish 11 types of instructions. Some important {\em types of instructions} are given in Overview \ref{Sigma_InstructionsFini}. 

\begin{overview}[$\sigma$-instructions and other instructions]\label{Sigma_InstructionsFini}

\hfill

\nopagebreak 

\noindent \fbox{\parbox{11.8cm}{

Computation instructions

\hspace{0.6cm} $(1)$ \quad $\ell : \, Z_j:= f_i^{m_i}(Z_{j_1},\ldots, Z_{j_{m_i}}) $

\hspace{0.6cm} $(2)$ \quad $\ell : \,Z_j:=c_i^0$

\vspace{0.1cm}

Copy instructions

 $^{\rm direct}(3)$ \quad $\ell : \,Z_{j}:=Z_{k}$

\vspace{0.1cm}

Branching instructions

\hspace{0.6cm} $(4)$ \quad {\sf $\ell : \,$ if $r_i^{k_i}(Z_{j_1},\ldots, Z_{j_{k_i}})$ then goto $\ell _1$ else goto $\ell _2$}\hspace*{0.5cm}

\vspace{0.1cm}

Stop instruction

\hspace{0.6cm} $(8)$ \quad$l : $ {\sf stop}}}
\end{overview}
We permit only one stop instruction in ${\cal P}$ and assume that it is of type (8) and that $l=\ell_{\cal P}$. The form of the instructions of types (1), (2), and (4) depends on the signature. All indices $j, k, j_1,j_2,\ldots $ are placeholders for positive integers, each $\ell$ is a placeholder for a label in $\{1,\ldots,l-1\}$, $\ell_1$ and $\ell_2$ are placeholders for labels in $\{1,\ldots,l\}$, and every $i$ stands for a positive integer that is less than or equal to $ n_1$, $n_2$, and $ n_3$, respectively. $ f_i^{m_i}$ is a symbol for an $m_i$-ary operation $ f_i:U_{\cal A}^{m_i}\to U_{\cal A}$, and so on. All allowed instructions can be substrings of a program and the dots $\ldots$ are not part of these substrings in the program. We write $f(Z_{j_1},\ldots, Z_{j_{m_i}})$ instead of $f((Z_{j_1},\ldots, Z_{j_{m_i}}))$ and we have $(Z_{j_1},\ldots, Z_{j_{1}})= Z_{j_1}$, $(Z_{j_1},\ldots, Z_{j_{2}})= (Z_{j_1},Z_{j_{2}})$, $(Z_{j_1},\ldots, Z_{j_{3}})= (Z_{j_1},Z_{j_{2}}, Z_{j_{3}})$, and so on.

The execution of a $\sigma$-program requires an interpretation of the constant, operation, and relation symbols and is possible by means of a structure ${\cal A}$ of signature $\sigma$. A machine using ${\cal A}$ is also called {\it machine over ${\cal A}$} or {\it ${\cal A}$-machine}. 
\begin{agree}[Restricted use of the identity]\label{AgreeIdent}An ${\cal A}$-machine can execute an equality test only if the identity $=$ restricted to $U_{\cal A}$ is a binary relation $r_i$ (with $i \leq n_3$) that belongs to ${\cal A}$ and contains a pair $(x_1,x_2)$ of individuals $x_1$ and $x_2$ if and only if $x_1=x_2$ holds in ${\cal A}$. It can be used for interpreting a condition such as $r_i^{2}(Z_{j_1}, Z_{j_2})$ that can be expressed less formally by $Z_{j_1} = Z_{j_2}$.
\end{agree}

\subsection{Finite machines and more} 

Let ${\cal A}$ be a structure of signature $\sigma$.
For any finite ${\cal A}$-machine ${\cal M}$, the number $j_{\cal M}$ of $Z$-registers is finite. Its program is a $\sigma$-program whose instructions can be of the types given in Overview \ref{Sigma_InstructionsFini} provided that the indices in the names of the $Z$-registers in these instructions satisfy the inequalities $j, k, j_1,j_2,\ldots \leq j_{\cal M}$.

\begin{overview}
[The registers of a finite ${\cal A}$-machine]\label{RegistFinDim}
\hfill

\nopagebreak 

\noindent \fbox{\parbox{11.8cm}{
\centering\vspace{0.1cm}
{\small
\begin{tabular}{c}

\noindent\noindent\vspace{0.2cm} \begin{tabular}{|c|}
\hline
\begin{tabular}{c|}$Z_1\,$\\\end{tabular}
\begin{tabular}{c|}$Z_2\,$\\\end{tabular}
\begin{tabular}{c|}$Z_3\,$\\\end{tabular}
\begin{tabular}{c|}$Z_4\,$\\\end{tabular}
\begin{tabular}{c|}$\ldots$\\\end{tabular}
\begin{tabular}{c}$\!\!\!Z_{j_{\cal M}}\!\!$\\\end{tabular}\\ 
\hline
\end{tabular}\hspace{0.5cm}\hfill{$Z$-registers (for {\em individuals} in $U_{\cal A}$)\,\,}\\
\noindent\noindent\vspace{0.1cm} \begin{tabular}{|c|}
\hline
\begin{tabular}{c}$B$\\\end{tabular}\\ 
\hline
\end{tabular}\hfill{ A register (for storing a {\em label}, an instruction counter)\,\,}\\
\end{tabular}}
}}
\end{overview}
The possible {\em change of the configurations} that are given by $(\ell,u_1,\ldots, u_{j_{\cal M}}) \in \{1,\ldots, \ell_{{\cal P}_{\cal M}}\}\times U_{\cal A} \times\cdots \times U_{\cal A}$ can be defined as in \cite{B92} and analogously to the change of the configurations of infinite machines defined below. 

To check whether $\sum_{i=1}^3 x_i=1$ holds for any numbers $x_1,x_2$, and $x_3$ in $\mathbb{R}^==_{\rm df }(\mathbb{R};1,0;+,-,\cdot;=)$, it is possible to use a finite $\mathbb{R}^=$-machine. Example \ref{Example1} shows how three numbers can be added.
If we want to check whether $\sum_{i=1}^{n} x_i=1$ holds for any tuple $(x_1,\ldots,x_n)$, we want to use an infinite dimensional machine such as given in Example \ref{VieleSummanden}. These machines are equipped with an infinite number of $Z$-registers and several index registers for storing indices such as $i$ and the length $n$. In the context of abstract computation over ${\cal A}=(U_{\cal A};\emptyset;; \in,=)$ where $U_{\cal A}$ is a set of sets, it is possible (as in Example \ref{Example2}) to check by a finite ${\cal A}$-machine whether a quantifier-free first-order formula $\phi(X_1,\ldots,X_n)$ of signature $(1;;2,2)$ becomes true in ${\cal A}$ when any sets $x_1,\ldots,x_n\in U_{\cal A}$ are assigned to $Z$-registers and we imagine that these registers stand for the first-order variables $X_1,\ldots,X_n$. If the truth values are to be determined for any assignment, it is possible to return codes, for instance, $\emptyset\in U_{\cal A}$ for $true$ and $(\emptyset,\emptyset)\in U_{\cal A}^2$ for $false$. This means that it is decidable by an ${\cal A}$-machine ${\cal M}_{\phi}$ for which assignments a given formula $\phi$ becomes true.
If we want to evaluate all formulas of a signature $\sigma$ and determine whether these formulas can be true in a structure ${\cal A}$ of signature $\sigma$ by using one single ${\cal A}$-machine, we need suitable codes for these formulas and we generally need more than a finite number of registers. One possibility is to take the codes of finite machines, each of which allows to evaluate one quantifier-free formula, as a part of an input of an infinite universal machine as defined in \cite{GASS20} for simulating the individual machines. The quantified formulas can be checked be non-deterministic universal oracle machines discussed later.

\subsection{Infinite machines and their configurations}\label{SecInfMa}

An infinite ${\cal A}$-machine is equipped with an infinite number of $Z$-registers $Z_1,Z_2$, $Z_3,\ldots$ for storing any arbitrary finite number of individuals or even an infinite sequence of individuals $u_1,u_2,\ldots \in U_{\cal A}$. We distinguish between memory and accumulator registers (which are also important components of a real RAM). For performing the instructions of types (1), (2), and (4) where the $Z$-registers to be used are uniquely determined, each ${\cal A}$-machine needs only a finite number $m_0$ of so-called {\em accumulator registers}. If the program of an infinite machine is a $\sigma$-program, then $m_0\geq \max\{2, m_1+1, \dots,m_{n_2}+1,k_1, \dots,k_{n_3}\}$ can be a necessary condition. Most of the $Z$-registers, $Z_{m+1},Z_{m+2},\ldots$ (for some $m\geq 0$), are so-called {\em memory registers} that form an extended memory-only store used for storing individuals. 
To be able to use all of these $Z$-registers during the computation process itself, the transfer between the memory registers and the accumulator registers can be necessary. It is possible due to the replacement of the instructions of type $^{\rm direct}(3)$ by instructions that allow to address registers indirectly. 
We follow the suggestion in \cite{BSS89} as far as possible and use as few types of instructions as possible. However, for calculating and evaluating addresses --- that are also called {\em indices} and that belong to the Peano structure ${\cal A}_{\bbbn}=_{\rm df}(\bbbn_+; 1;succ;=)$ defined by $\bbbn_+=_{\rm df}\bbbn\setminus\{0\}$ and $succ(n)=n+1$ ---, we allow any finite number of index registers. We do this for the following reasons. On the one hand, we want to save the hassle of managing complicated information exchange between index registers and the $Z$-registers and the use of operations of ${\cal A}$ for computing addresses. On the other hand, we want to be able to take structures ${\cal A}$ into account that do not contain a successor function such as $succ$ and that cannot be expanded to a structure ${\cal A}'$ with $U_{\cal A}=U_{{\cal A}'}$ by adding a function computable over ${\cal A}$ such that ${\cal A}_{\bbbn}$ is isomorphic to a substructure of ${\cal A}'$.
\begin{overview}
[The registers of an infinite ${\cal A}$-machine ${\cal M}$]\label{RegistInfinDim}
\hfill

\nopagebreak 

\noindent \fbox{\parbox{11.8cm}{
\centering\vspace{0.1cm}
{\small
\begin{tabular}{c}
\noindent\vspace{0.2cm} \begin{tabular}{|c}
\hline
\begin{tabular}{c|}$Z_1$\\\end{tabular}
\begin{tabular}{c|}$Z_2$\\\end{tabular}
\begin{tabular}{c|}$Z_3$\\\end{tabular}
\begin{tabular}{c|}$Z_4$\\\end{tabular}
\begin{tabular}{c|}$Z_5$\\\end{tabular}
\begin{tabular}{c}$\ldots$\\\end{tabular}\\
\hline
\end{tabular} 
\hspace{0.5cm}\hfill{$Z$-registers (for {\em individuals} in $U_{\cal A}$)}\\
\noindent\noindent\vspace{0.2cm} \begin{tabular}{|c|}
\hline
\begin{tabular}{c|}$I_1\!$\\\end{tabular}
\begin{tabular}{c|}$I_2$\\\end{tabular}
\begin{tabular}{c|}$I_3$\\\end{tabular}
\begin{tabular}{c|}$I_4$\\\end{tabular}
\begin{tabular}{c|}$\ldots$\\\end{tabular}
\begin{tabular}{c}$\!\!I_{k_{\cal M}}\!\!$\\\end{tabular}\\ 
\hline
\end{tabular}\hspace{1cm}\hfill{Index registers (for {\em indices} in $\mathbb{N}_+$)}\\
\noindent\noindent\vspace{0.1cm} \begin{tabular}{|c|}
\hline
\begin{tabular}{c}$B$\\\end{tabular}\\ 
\hline
\end{tabular}\hfill{\qquad A register (for storing a {\em label}, an instruction counter)}\\
\end{tabular}}
}}
\end{overview} 

Let ${\cal M}$ be an infinite ${\cal A}$-machine equipped with $k_{\cal M}$ {\em index registers} $I_1,\ldots, I_{k_{\cal M}}$ for some $k_{\cal M}\geq 1$.
The execution of instructions of types (3) and (5) to (7) given in Overview \ref{Index_Instructions} by ${\cal M}$ is possible only if the so-called {\em names of index registers} in the considered instructions belong to $I_1,\ldots, I_{k_{\cal M}}$. While $^{\rm direct}$(3) describes copying by using {\em direct addressing}, the instructions of type (3) allow to read the content of a memory register and to write a value in a memory register by using {\em indirect addressing}. For any program ${\cal P}$, let $k_{\cal P}$ be the smallest integer $k_0$ with $j\leq k_0$ for any index $j$ that occurs in a name of an index register in ${\cal P}$.
\begin{overview}
[Instructions for the transport of objects]\label{Index_Instructions}

\hfill

\nopagebreak 

\noindent \fbox{\parbox{ 11.8cm}{

Copy instructions

\hspace{0.6cm} $(3)$ \quad $\ell : \,Z_{I_j}:=Z_{I_k}$

\vspace{0.1cm}

Index instructions

\hspace{0.6cm} $(5)$ \quad {\sf $\ell : \,$ if $I_{j}=I_k$ then goto $\ell _1$ else goto $\ell _2$} 

\hspace{0.6cm} $(6)$ \quad $\ell : \,I_j:=1$

\hspace{0.6cm} $(7)$ \quad $\ell : \,I_j:=I_j+1$
}}
\end{overview}
The execution of a program ${\cal P}$ containing instructions of types (1) to (8) is possible by an infinite ${\cal A}$-machine ${\cal M}$ only if the constant, function and relation symbols can be interpreted by means of ${\cal A}$ and $k_{\cal M}\geq k_{\cal P}$ holds.

Any infinite ${\cal A}$-machine ${\cal M}$ is an abstract machine. The {\em execution of an instruction} in ${\cal P}_{\cal M}$ by ${\cal M}$ can be described purely mathematically by single-valued or multi-valued functions whose application can lead to a new overall state of the machine that is called a {\em configuration} of ${\cal M}$. At any point of time, the configuration of ${\cal M}$ results from the values stored in the registers of ${\cal M}$. Let $U^\omega=_{\rm df}\{(u_1,u_2,\ldots)\mid (\exists u:\bbbn_+ \to U)(\forall i\in \bbbn_+)(u(i)=u_i)\}$ for any set $U$. Any configuration is given by a sequence
 $(\ell, \nu_1,\ldots,\nu_{k_{\cal M}}, u_1,u_2,\ldots)$ denoted shortly by $(\ell\,.\,\vec \nu\,.\,\bar u)$ with $\ell \in \{1,\ldots, \ell_{{\cal P}_{\cal M}}\}$, $\vec \nu=( \nu_1,\ldots,\nu_{k_{\cal M}})\in \mathbb{N}_+^{k_{\cal M}}$, and $\bar u=( u_1,u_2,\ldots)\in U_{\cal A}^\omega$. The label $\ell$ is the {\em content $c(B)$ of $B$}. For each $i\leq k_{\cal M}$, the positive integer $\nu_i$ is the {\em content $c(I_i)$ of $I_i$}. For any $i\geq 1$, the individual $u_i$ is the {\em content $c(Z_i)$ of $Z_i$}. $\{ (\ell\,.\,\vec \nu\,.\, \bar u)\mid 1\leq \ell\leq \ell_{{\cal P}_{\cal M}} \,\,\&\,\, \vec \nu\in (\mathbb{N}_+)^{k_{\cal M}} \,\,\&\,\, \bar u\in U_{\cal A} ^\omega\}$ is the {\em space ${\sf S}_{\cal M}$ of all possible configurations of ${\cal M}$}. 
A so-called {\em input procedure} is used to determine the {\em initial configuration} of ${\cal M}$.
The transformation of a configuration $(\ell\,.\,\vec \nu\,.\,\bar u)$ into the next configuration is {\em one computation step of ${\cal M}$}. It is determined by applying
 the instruction in ${\cal P}_{\cal M}$ whose label is $\ell$ and completely defined by a binary relation $\to_{\cal M}\, \subseteq {\sf S}_{\cal M}\times {\sf S}_{\cal M}$ resulting from single-valued or multi-valued functions in a system ${\cal F}_{\cal M}$. For examples, see Figure \ref{InstructionsChangeConfInfinite}. The full {\em transition system ${\cal S}_{\cal M}=({\sf S}_{\cal M},\to_{\cal M} )$} is defined in \cite{GASS20}. It allows to define each possible sequence of transformations of configurations of ${\cal M}$ step by step. If ${\cal M}$ can execute only instructions of types (1) to (8), $\to_{\cal M}$ is a {\em total function} of the form $\to_{\cal M}: {\sf S}_{\cal M}\to {\sf S}_{\cal M}$. According to this, there is a {\em partial function} $(\to_{\cal M})_{{\rm Stop}_{\cal M}}: \, \subseteq {\sf S}_{\cal M}\to {\sf S}_{\cal M}$ assigning to any initial configuration $con_1$ of ${\cal M}$ a
 stop configuration
$(\ell_{\cal P}\,.\,\vec \nu\,.\,\bar u)$
 if and only if the execution of instructions of ${\cal P}_{\cal M}$ by $\cal M$, that starts with $con_1$ and that is determined step by step by $\to_{\cal M}$, results in this stop configuration. An example is given in \,Example \ref{PositivOrakel}. Formally, we can describe an infinite dimensional ${\cal A}$-machine ${\cal M}$ as follows (cf. \cite{GASS11, GASS20}).

\begin{overview}[A deterministic machine]\label{Components}
\hfill

\nopagebreak 
\noindent \fbox{\parbox{11.8cm}{
A deterministic machine
$(U ^\omega,(\mathbb{N}_+)^k, {\cal L}, {\cal P},{\cal B} , {\rm In}, {\rm Out})$ over ${\cal A}$ consists of
\begin{itemize}\itemsep0pt
\item a space of memory states, $U ^\omega$, 
\item a space of addresses (or indices), $(\mathbb{N}_+)^k$, with $k\geq 1$, 
\item a set of labels, ${\cal L}=\{1,\ldots, l\}$, with $l\geq 1$,
\item a program ${\cal P}\in {\sf P}_{\sigma}$ with $\ell_{\cal P}= l$ and $k_{\cal P}\leq k$,
\item a structure ${\cal B}=(U;c_{\alpha_1},\ldots,c_{\alpha_{n_1}};f_1,\ldots, f_{n_2}; r_1,\ldots,r_{n_3})$ of signature $\sigma$,
\item an input function ${\rm In}: {\sf I}\to\{ (\vec\nu\,.\,\bar u)\mid(\vec\nu,\bar u)\in(\mathbb{N}_+)^k\times U^{\omega}\}$,
\item an output function ${\rm Out}: \{ (\vec\nu\,.\,\bar u)\mid (\vec\nu,\bar u)\in(\mathbb{N}_+)^k\times U^{\omega}\}\to {\sf O}$.
\end{itemize}}}
\end{overview}
A tuple 
$(U_{\cal } ^\omega,(\mathbb{N}_+)^k, {\cal L}, {\cal P},{\cal B} , {\rm In}, {\rm Out})$ whose components are determined as in Overview \ref{Components} is an {\em infinite dimensional machine over ${\cal A}$} (or {\em ${\cal A}$-machine}) if ${\cal B}$ is a reduct of ${\cal A}$ (which means that $U$ stands for $U_{\cal A} $). Let ${\cal M}$ be given by ${\cal M}=(U_{\cal A} ^\omega,(\mathbb{N}_+)^k, {\cal L}, {\cal P},{\cal B} , {\rm In}, {\rm Out})$, let $k_{\cal M}=k$, ${\cal L}_{\cal M}={\cal L}$, and ${\cal P}_{\cal M}={\cal P}$. The constants $c_{\alpha_1},\ldots,c_{\alpha_{n_1}}$ of ${\cal B}$ can be the {\em machine constants of ${\cal M}$}. 
Let ${\rm In}_{\cal M}$ and ${\rm Out}_{\cal M}$ be the input function ${\rm In}$ and the output function ${\rm Out}$, respectively, of ${\cal M}$. Let ${\sf I}_{\cal M}$ be the {\em input space of ${\cal M}$} given by ${\sf I}_{\cal M}={\sf I}$, and ${\sf O}_{\cal M}$ be the {\em output space of ${\cal M}$} satisfying ${\sf O}= {\sf O}_{\cal M}$. Here and in the following, $(\vec \nu\,.\,\bar u)$ stands for $(\nu_1, \ldots,\nu_{k_{\cal M}},u_1,u_2,\ldots)$ and $(\vec \nu,\bar u)$ stands for $((\nu_1, \ldots,\nu_{k_{\cal M}}),(u_1,u_2,\ldots))$, and so on. The {\em instructions in ${\cal P}_{\cal M}$} can (initially \footnote{Later, we also allow instructions of other types.}) be of types (1) to (8). The functions collected in a {\em system ${\cal F}_{\cal M}$} for defining the relation $\to_{\cal M}$ and for transforming the configurations of ${\cal M}$ are derived from ${\cal P}_{\cal M}$ and its interpretation. Their definitions are given in Overviews \ref{Oper1}, \ref{copy}, \ref{Test1}, and \ref{IndOper}. For more details, see \cite{GASS20}.

\begin{overview}[Elementary operations in ${\cal F}_{\cal M}$ for changing the state]\label{Oper1}
\hfill

\nopagebreak 
\noindent \fbox{\parbox{11.8cm}{
\begin{itemize}\itemsep0pt
\item
$F_\ell:U_{\cal A} ^{\omega}\to U_{\cal A} ^{\omega}$ \hfill for instructions of type $(1)$
\parskip -0.5mm 
\begin{tabbing} 
$C_{\ell }(\vec \nu,\bar u)\quad$\=$=( u_1,\ldots,u_{j-1},f_i(u_{j_1},\ldots, u_{j_{m_i}}),u_{j+1}, u_{j+2}, \ldots)$ \= \kill
$ F_{\ell }(\bar u)$\>$=( u_1,\ldots,u_{j-1},f_i(u_{j_1},\ldots, u_{j_{m_i}}),u_{j+1}, u_{j+2}, \ldots)$
\end{tabbing} 
\item
$F_\ell:U_{\cal A} ^{\omega}\to U_{\cal A} ^{\omega}$ \hfill for instructions of type $(2)$
\parskip -0.5mm 
\begin{tabbing} 
$C_{\ell }(\vec \nu,\bar u)\quad$\=$=( u_1,\ldots,u_{j-1},f_i(u_{j_1},\ldots, u_{j_{m_i}}),u_{j+1}, u_{j+2}, \ldots)$ \= \kill
$ F_{\ell }(\bar u)$\>$=( u_1,\ldots,u_{j-1},c_{\alpha_i},u_{j+1}, u_{j+2}, \ldots) $
\end{tabbing} 
\end{itemize}}}
\end{overview}

\begin{overview}[Copy functions in ${\cal F}_{\cal M}$ for changing the state]\label{copy}
\hfill

\nopagebreak 
\noindent \fbox{\parbox{11.8cm}{
\begin{itemize}\itemsep0pt
\item $C_\ell: (\mathbb{N}_+)^{k_ {\cal M}} \times U_{\cal A} ^{\omega}\to U_{\cal A} ^{\omega}$ \hfill for instructions of type $(3)$
\parskip -0.5mm 
\begin{tabbing} $C_{\ell }(\vec \nu,\bar u)\quad$\=$=( u_1,\ldots,u_{j-1},f_i(u_{j_1},\ldots, u_{j_{m_i}}),u_{j+1}, u_{j+2}, \ldots)$ \= \kill
$C_{\ell }(\vec \nu,\bar u)$\>$=(u_1,\ldots,u_{\nu_j-1},u_{\nu_k},u_{\nu_j+1}, u_{\nu_j+2}, \ldots) $
\end{tabbing} 
\end{itemize}}}
\end{overview}

\begin{overview}[Test functions in ${\cal F}_{\cal M}$ for changing the label] \label{Test1}
\hfill

\nopagebreak 
\noindent \fbox{\parbox{11.8cm}{
\begin{itemize}\itemsep0pt
\item 
$T_\ell: U_{\cal A} ^{\omega} \to {\cal L}_ {\cal M}$ \hfill for instructions of type $(4)$
\parskip -0.5mm 
\begin{tabbing} $C_{\ell }(\vec \nu,\bar u)\quad$\=$=( u_1,\ldots,u_{j-1},f_i(u_{j_1},\ldots, u_{j_{m_i}}),u_{j+1}, \ldots)$ \= \kill
$T_{\ell }(\bar u)$\>$=\left\{\begin{array}{ll}
\ell_1\quad\mbox{ if $ (u_{j_1},\ldots, u_{j_{k_i}})\in r_i$}\\
\ell_2\quad\mbox{ if $ (u_{j_1},\ldots, u_{j_{k_i}})\not\in r_i$}\end{array}\right.$
\end{tabbing} 
\item
$T_\ell: (\mathbb{N}_+)^{k_ {\cal M}} \to {\cal L}_ {\cal M}$ \hfill for instructions of type $(5)$
\parskip -0.5mm 
\begin{tabbing} $C_{\ell }(\vec \nu,\bar u)\quad$\=$=( u_1,\ldots,u_{j-1},f_i(u_{j_1},\ldots, u_{j_{m_i}}),u_{j+1}, \ldots)$ \= \kill
$T_{\ell }(\vec \nu)$\>$=\left\{\begin{array}{ll}
\ell_1\quad \mbox{ if $ \nu_j=\nu_k$}\\
\ell_2\quad \mbox{ if $ \nu_j\not=\nu_k$}\end{array}\right.$ 
\end{tabbing} 
\end{itemize}}}
\end{overview}

\begin{overview}[Auxiliary functions in ${\cal F}_{\cal M}$ for changing an index]\label{IndOper}
\hfill

\nopagebreak 
\noindent \fbox{\parbox{11.8cm}{
\begin{itemize}\itemsep0pt
 \item {\em } $H_\ell: (\mathbb{N}_+)^{k_ {\cal M}} \to (\mathbb{N}_+)^{k_ {\cal M}} $ \hfill for instructions of type $(6)$ 
\parskip -0.5mm 
\begin{tabbing} $C_{\ell }(\vec \nu,\bar u)\quad$\=$=( u_1,\ldots,u_{j-1},f_i(u_{j_1},\ldots, u_{j_{m_i}}),u_{j+1}, \ldots)$ \= \kill
$H_{\ell }(\vec \nu)$\>$= ( \nu_1,\ldots,\nu_{j-1},1,\nu_{j+1}, \ldots,\nu_{k_{\cal M}})$ 
\end{tabbing} 
\item $H_\ell: (\mathbb{N}_+)^{k_ {\cal M}} \to (\mathbb{N}_+)^{k_ {\cal M}} $ \hfill for instructions of type $(7)$
\parskip -0.5mm 
\begin{tabbing} $C_{\ell }(\vec \nu,\bar u)\quad$\=$=( u_1,\ldots,u_{j-1},f_i(u_{j_1},\ldots, u_{j_{m_i}}),u_{j+1}, \ldots)$ \= \kill
$H_{\ell }(\vec \nu)$\>$= ( \nu_1,\ldots,\nu_{j-1},\nu_j+1,\nu_{j+1}, \ldots,\nu_{k_{\cal M}})$ 
\end{tabbing} 
\end{itemize}}}
\end{overview}
Consequently, each instruction causes the {\em change of a configuration of ${\cal M}$} resulting from applying a {\em transition rule} given by $\to_{\cal M}$ (as in Figure \ref{InstructionsChangeConfInfinite}). A pair $((\ell_{{\cal P}_{\cal M}}\,.\,\vec \nu\,.\,\bar u),(\ell\,.\,\vec \mu\,.\,\bar z))\in ({\sf S}_{\cal M})^2$ belongs to $\to_{\cal M}$ only if $\ell=\ell_{{\cal P}_{\cal M}}$, $\vec \nu=\vec \mu$, and $\bar u=\bar z$. If $\ell_1$ is the label of an instruction of types (1) to (3) or (6) or (7), then $((\ell_1\,.\,\vec \nu\,.\,\bar u), (\ell_2\,.\,\vec \mu\,.\,\bar z))$ can belong to $\to_{\cal M}$ only if $\ell_2=\ell_1+1$. Any pair $((1\,.\,\vec \nu\,.\,\bar u), (\ell_2\,.\,\vec \mu\,.\,\bar z))$ in the transitive closure of $\to_{\cal M}$ belongs to $(\to_{\cal M})_{{\rm Stop}_{\cal M}}$ if and only if $\ell_2=\ell_{{\cal P}_{\cal M}}$. The {\em computation of ${\rm Res}_{\cal M}:\, \subseteq\! {\sf I}_{\cal M}\to {\sf O}_{\cal M}$ by ${\cal M}$ } is determined by a {\sf computational system ${\cal R}_{\cal M}$} consisting of ${\cal S}_{\cal M}$, an {\sf input} and an {\sf output procedure}, ${\rm Input}_{\cal M}:{\sf I}_{\cal M}\to {\sf S}_{\cal M}$ and ${\rm Output}_{\cal M}:{\sf S}_{\cal M} \to {\sf O}_{\cal M}$, and a {\sf stop criterion}. We have ${\rm {\rm Res}_{\cal M}}({\sf i})={\rm Output}_{\cal M}\circ (\to_{\cal M})_{{\rm Stop}_{\cal M}} \circ {\rm Input}_{\cal M}({\sf i})$ for all $\,{\sf i}\in {\sf I}_{\cal M}$. 
\begin{figure}
\centering
{\begin{tabular}{|l|} \hline\\
$((\ell \,.\,\vec \nu\,.\,\bar u) ,(\ell +1 \,.\,\vec \nu\,.\,\underbrace {(u_1,\ldots,u_{\nu_j-1},u_{\nu_k},u_{\nu_j+1},u_{\nu_j+2}, \ldots)}_{ C_{\ell }(\vec \nu,\bar u) }\,))\in \,\to_{\cal M} $\\ \hspace{7.5cm} caused by $\ell: \,Z_{I_{j}}:=Z_{I_k}$ \\\\
$((\ell \,.\,\vec \nu\,.\,\bar u),(\ell +1\,.\, \underbrace {( \nu_1,\ldots,\nu_{j-1},\nu_j+1,\nu_{j+1}, \ldots,\nu_{k_{\cal M}})}_{ H_{\ell}(\vec \nu)}\,.\,\bar u))\in \,\to_{\cal M} $ \\\hfill caused by $\ell: \,I_{j}:=I_j+1$ \\\hline \end{tabular}}
\caption{Examples for pairs $(con_s,con_{s+1})$ of configurations in $\to_{\cal M}$ } \label{InstructionsChangeConfInfinite}
\end{figure}
\begin{consequ}[A side effect]\label{SideEffect} Note, that the instructions of types (5) to (7) are primarily intended to enable copying. However, it is also possible to simulate machines over ${\cal A}_{\mathbb{N}}$ by executing instructions of these types. The code of the values stored in two index registers of an ${\cal A}$-machine can be stored in one single register by using the restricted Cantor pairing function $cantor:\mathbb{N}_+^2\to \mathbb{N}_+$ denoted also by $ca$ and defined by $ca(\mu_1,\mu_2)=_{\rm df} \frac12((\mu_1+\mu_2)^2 +3\mu_2+\mu_1)$.\linebreak We know that the values $cantor(\mu_1,\mu_2)$ as well as $ca_1(ca(\mu_1,\mu_2))=_{\rm df}\mu_1$ and $ca_2(ca(\mu_1,\mu_2))=_{\rm df}\mu_2$ computable over ${\cal A}_{\mathbb{N}}$ are also computable with the help of some index registers and suitable index instructions of types (5) to (7).
\end{consequ}

\subsection{Computation and semi-decidability by BSS RAMs}\label{Det_Non_Det_BSS_RAM}

A known model of computation is the abstract Turing machine (TM) with a tape consisting of an infinite number of cells for storing at least two symbols. Turing machines that use only the symbols 0 and 1, and a blank symbol can be {\sf simulated} by ${\cal A}_0$-machines with ${\cal A}_0=_{\rm df}(\{0,1\};0,1;;=)$ as described in \cite[pp.\,589--590]{GASS20}. A Turing-computable function $f:\,\subseteq\mathbb{N}_+\to \mathbb{N}_+$ can be computed by an ${\cal A}_0$-machine by using an {\sf input procedure} that makes it possible to assign the symbols (the digits) of a (binary) string encoding any positive integer to some of the $Z$-registers and an {\sf output procedure} that, after the simulation of a TM, decodes the symbols stored in the $Z$-registers into the corresponding natural number if it is possible. (For more, see also Example \ref{Example0TM}).

 Roughly speaking, a function $f:\,\subseteq\mathbb{R}\to \mathbb{R}$ is computable by a Type-2 TM $M$ if, for any $x\in \mathbb{R}$, the input procedure ${\rm input}_M$ provides suitable infinite sequences in $\{0,1\}^\omega$ as codes for $x$ that are processed by a program to obtain sequences in $\{0,1\}^\omega$ that could be decoded by means of an output procedure ${\rm output}_M$ into a real number $f(x)$. For more, see \cite{WEIHRAUCH}.

\begin{overview}[Type-2 computation of $f:\,\subseteq\mathbb{R}\to \mathbb{R}$ over ${\cal A}_0$]

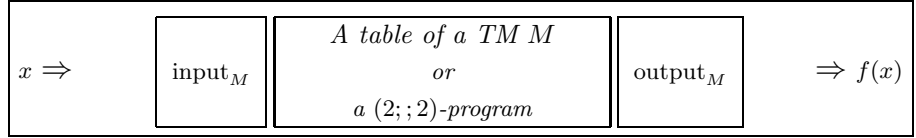
\hfill

\nopagebreak 
\noindent \fbox{\parbox{11.8cm}{
\centering\vspace{0.1cm}
\small
$ x$ {\large $\Rightarrow$ \qquad} 
{\small
\begin{tabular}{|c|}\hline\\${\rm input}_ M$\\\\\hline\end{tabular}\hspace{0.1cm}
$\!\!$\begin{tabular}{|c|}\hline
\hspace*{0.4cm} {\normalsize A table of a TM $M$} \hspace*{0.4cm}\\or\\a $(2;;2) $-program \\\hline
\end{tabular}\hspace{0.1cm}\begin{tabular}{|c|}\hline\\
${\rm output}_ M$\\\\\hline\end{tabular}} {\large \qquad $\Rightarrow$} $f( x)$
}}
\end{overview}

In contrast to type-2 machines over ${\cal A}_ 0$ or over $(U;\emptyset;\cup;=)$ whose outputs generally result from decoding sequences in $\{0,1\}^\omega$ or $U^\omega$ that can, for instance, be the code of a real number $r\in \mathbb{R}$ or a set $\bigcup_{i=1}^\infty [0,\frac{i-1}{i}[$ (with $[0,\frac{i-1}{i}[\in U$), a BSS RAM ${\cal M}$ provides outputs by {\sf applying its output procedure only if} ${\cal M}$ halts.
\begin{overview}[The BSS-RAM model]
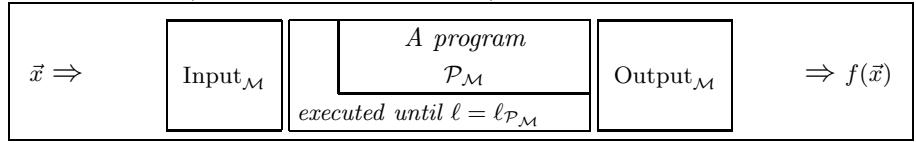
\label{BSSRAMModell}

\hfill

\nopagebreak 
\noindent \fbox{\parbox{11.8cm}{
\centering\vspace{0.1cm}
{\small
$\vec x$ {\large $\Rightarrow$ \qquad} 
{\small
\begin{tabular}{|c|}\hline\\${\rm Input}_{\cal M}$\\\\\hline\end{tabular}\hspace{0.1cm}\begin{tabular}{|r|}\hline \,\,\,
\begin{tabular}{|c}
\hspace*{0.54cm} {\normalsize A program} \hspace*{0.54cm}\\${\cal P}_{\cal M}$\\\hline
\end{tabular}$\!\!\!\!$\\$\!\!$executed until $\ell=\ell_{{\cal P}_{\cal M}}$ \quad \qquad \\\hline
\end{tabular}\hspace{0.1cm}\begin{tabular}{|c|}\hline\\
${\rm Output}_{\cal M}$\\\\\hline\end{tabular}} {\large \qquad $\Rightarrow$} $f(\vec x)$
}}}
\end{overview}

Let ${ U}^{\infty}=_{\rm df} \bigcup_{n\geq 1} { U}^n$ for any set $U$. For any $\vec x\in U^\infty$, let $n$ denote the respective length of $\vec x$ such that we have $length (\vec x)=n$, $\vec x=(x_1,\ldots, x_n)$, and $n\in \bbbn_+$. An ${\cal A}$-machine ${\cal M}$ whose {\sf input space} ${\sf I}_{\cal M}$ and whose {\sf output space} ${\sf O}_{\cal M}$ are ${ U}^{\infty}_{\cal A}$ is called a {\em deterministic BSS RAM over ${\cal A}$} if its {\sf input procedure} ${\rm Input}_{\cal M}: U^{\infty}_{\cal A} \to {\sf S}_{\cal M}$ is the function given by
\begin{equation}\tag{I1}\label{Inp1}{\rm Input}_{\cal M}(x_1,\ldots, x_n)= (1\,.\,(n, 1,\ldots, 1)\,.\,(x_1, \ldots,x_n,x_n,x_n, \ldots )),
\end{equation}
for all $\vec x\in U_{\cal A}^\infty$, where $(n, 1,\ldots, 1)$ stands for a tuple in $(\mathbb{N}_+)^{k_{\cal M}}$, and its {\sf output procedure} ${\rm Output}_{\cal M}: {\sf S}_{\cal M}\to U^{\infty}_{\cal A} $ is the function given by 
\begin{equation}\tag{O}
\label{Out}{\rm Output}_{\cal M}(\ell\,.\,\vec \nu\,.\,\bar u)= (u_1,\ldots, u_{\nu_1})
\end{equation}
for any configuration $(\ell\,.\,(\nu_1,\ldots,\nu_{k_{\cal M}})\,.\,(u_1,u_2,u_3,\ldots))\in{\sf S}_{\cal M} $. More precisely, these procedures are derived from the input and the output function as follows.
Let ${\rm Input}_{\cal M}: {\sf I}_{\cal M}\to {\sf S}_{\cal M}$ be defined by ${\rm Input}_{\cal M}({\sf i})= (1\,.\,{\rm In}_{\cal M}({\sf i}))$ for all ${\sf i}\in {\sf I}_{\cal M}$ such that ${\rm In}_{\cal M}$ is a unique map given by ${\rm In}_{\cal M}(\vec x)=((n, 1,\ldots, 1)\,.\, \vec x\,.\,(x_n,x_n, \ldots )) $ for all $\vec x\in {\sf I}_{\cal M}$. Let ${\rm Output}_{\cal M}: {\sf S}_{\cal M}\to {\sf O}_{\cal M}$ be defined by ${\rm Output}_{\cal M}(\ell\,.\,\vec \nu\,.\, \bar u)={\rm Out}_{\cal M}(\vec \nu\,.\, \bar u)$ and ${\rm Out}_{\cal M}(\vec \nu\,.\, \bar u) = (u_1,\ldots, u_{\nu_1})$ for all $(\ell\,.\,\vec \nu\,.\, \bar u)\in {\sf S}_{\cal M}$. 

Let ${\sf M}_{\cal A}$ be the class of {\sf all deterministic BSS RAMs} over ${\cal A}$ that can only execute instructions of {\sf types $(1)$ to $(8)$}. 

\begin{definition}[Computable function]$\!$For every BSS RAM ${\cal M}$ in ${\sf M}_{\cal A}$, 
 the partial function ${\rm Res}_{\cal M}:\,\subseteq\!U_{\cal A}^\infty \to U_{\cal A}^\infty $ defined by \begin{equation}\tag{R1}\label{Res1}{\rm {\rm Res}_{\cal M}}(\vec x)={\rm Output}_{\cal M}( (\to_{\cal M})_{{\rm Stop}_{\cal M}} ( {\rm Input}_{\cal M}({\vec x})))
\end{equation} 
for all $\vec x \in U_{\cal A}^\infty $ is called {\em the result function of $ {\cal M}$}. A partial function $f:\,\subseteq\! U_{\cal A}^\infty \to U_{\cal A}^\infty $ is {\em computable by a BSS RAM over ${\cal A}$} if it is the result function of a BSS RAM in ${\sf M}_{\cal A}$.
\end{definition} 
For ${\cal M}\in {\sf M}_{\cal A}$ and any $\vec x\in U_{\cal A}^\infty$, ${\cal M} ( \vec x)\downarrow $ means that ${\rm {\rm Res}_{\cal M}}(\vec x)$ is defined. In such a case, we say that ${\cal M}$ {\em halts} on {\sf input} $\vec x\in U_{\cal A}^\infty$ and {\em accepts} $\vec x$. For any {\sf output} given by ${\rm Res}_{\cal M}(\vec x)$, we have ${\rm Res}_{\cal M}(\vec x)\in U_{\cal A}^\infty$. Let $H_{\cal M}$ be the {\em halting set of ${\cal M }$} that is given by $H_{\cal M}=\{ \vec x \in U_{\cal A}^{\infty} \mid {\cal M} ( \vec x)\downarrow \}$. A subset $P\subseteq U_{\cal A}^\infty$ is called {\em semi-decidable} ({\em by a machine in ${\sf M}_{\cal A}$}) ({\em over ${\cal A}$}) if it is the halting set of a BSS RAM in $ {\sf M}_{\cal A}$. 
 Let ${\rm SDEC}_{\cal A}$ be the class of {\sf all problems $P\subseteq U_{\cal A}^\infty $ that are semi-decidable} by a BSS RAM in ${\sf M}_{\cal A}$ and let ${\rm DEC}_{\cal A}$ be the class of {\sf all decidable problems} $P\subseteq U_{\cal A}^\infty $ for which both the set $P$ and its complement $U_{\cal A}^\infty\setminus P$ are semi-decidable by two ${\cal A}$-machines in ${\sf M}_{\cal A}$.

Consequently, when ${\cal A}$ contains at least one constant $c_1, $ $P$ is in ${\rm SDEC}_{\cal A}$ if and only if its {\sf partial characteristic function} $\bar \chi_P: \, \subseteq U_{\cal A}^\infty\to \{c_1\}$ --- that is defined by $\bar \chi_P(\vec x)=c_1$ for all $\vec x \in P$ and $\bar \chi_P(\vec x)\uparrow$ (which means that $\bar \chi_P(\vec x)$ is not defined) for all $\vec x \in U_{\cal A}^\infty \setminus P$ --- is computable over ${\cal A}$. If ${\cal A}$ contains two constants $c_1$ and $c_2$ and $P\in {\rm DEC}_{\cal A}$ holds, then we can prove that its {\sf characteristic function} $\chi_P:U_{\cal A}^\infty\to \{c_1,c_2\}$ --- that is defined by $\chi_P(\vec x)=c_1$ for all $\vec x \in P$ and $\chi_P(\vec x)=c_2$ for all $\vec x \in U_{\cal A}^\infty \setminus P$ --- is computable over ${\cal A}$ and we can say that each $\vec x\in P$ is {\em accepted} and each $\vec x\in U_{\cal A}^\infty \setminus P$ is {\em rejected}. If ${\cal A}$ contains the identity or $\{c_1\}$ and $\{c_2\}$ are semi-decidable, then the computability of $\chi_P$ over ${\cal A}$ implies $P\in {\rm DEC}_{\cal A}$. We have, for example, $\{\vec x \in \mathbb{R}^\infty\mid \sum_{i=1}^{n} x_i=1\} \in {\rm SDEC}_{\mathbb{R}^=}$, and Example \ref{VieleSummanden} shows the decidability of $P=\{\vec x \in \mathbb{C}^\infty\mid \sum_{i=1}^{n} x_i=1\}$ over $\mathbb{C}=_{\rm df}(\mathbb{C};1,0,i;+,-,\cdot;=)$ by computing $\chi_P: \mathbb{C}^\infty\to \{1,0\}$. It is known that $\mathbb{R}_{\geq 0}=_{\rm df}\{x \in \mathbb{R} \mid x\geq 0\}$ is not semi-decidable over $\mathbb{R}^=$. Example \ref{PositivRaten} shows how $\mathbb{R}_{\geq 0}$ can be semi-decided by using a non-deterministic BSS RAM over $\mathbb{R}^=$.

\subsection{Non-deterministic BSS RAMs}
The difference between the deterministic BSS RAMs and the non-deterministic BSS RAMs that we now want to define results solely and exclusively from replacing the single-valued input function by a new multi-valued input function. 
\begin{overview}[An input function for a non-deterministic machine]\label{NonComponents}
\hfill

\nopagebreak 
\noindent \fbox{\parbox{11.8cm}{

A non-deterministic machine
$(U ^\omega,(\mathbb{N}_+)^k, {\cal L}, {\cal P},{\cal B} , {\rm In}, {\rm Out})$ can contain, e.g.,
\begin{itemize}\itemsep0pt
\item a multi-valued input function ${\rm In}: {\sf I}\genfrac{}{}{0pt}{2}{\longrightarrow} {\longrightarrow}\{ (\vec\nu\,.\,\bar u)\mid(\vec\nu,\bar u)\in(\mathbb{N}_+)^k\times U^{\omega}\}$

such that ${\rm In}\subseteq {\sf I}\times \{ (\vec\nu\,.\,\bar u)\mid(\vec\nu,\bar u)\in(\mathbb{N}_+)^k\times U^{\omega}\}$.
\end{itemize}}}
\end{overview}
Such an input function and the corresponding input procedure describe a guessing process. For any non-deterministic BSS RAM ${\cal M}$ over ${\cal A}$ considered here, the input function ${\rm In}_{\cal M}: U_{\cal A}^\infty \genfrac{}{}{0pt}{2}{\longrightarrow} {\longrightarrow} (\{ (\vec\nu\,.\,\bar u)\mid(\vec\nu,\bar u)\in(\mathbb{N}_+)^k\times U^{\omega}\}$ is given by \[{\rm In}_{\cal M}=\{(\vec x,(\,\underbrace{(n, 1,\ldots, 1)}_{\in(\mathbb{N}_+)^{k_{\cal M}}}\,.\,\underbrace{\vec x\,.\,\vec y\,.\,(x_n,x_n, \ldots) }_{\in U_{\cal A}^{\omega}}\,))\mid \vec x\in U_{\cal A}^\infty\,\,\&\,\,\vec y \in U_{\cal A}^\infty\}\] where $\vec x=(x_1,\ldots,x_n)$ and $(\vec x\,.\,\vec y)=(x_1,\ldots,x_n,y_1,\ldots,y_m)$ hold for $\vec x\in U_{\cal A}^n$ and $ \vec y\in U_{\cal A}^m$, and so on. Formally, the input procedure ${\rm Input}_{\cal M}:U_{\cal A}^\infty \genfrac{}{}{0pt}{2}{\longrightarrow} {\longrightarrow} {\sf S}_{\cal M} $ is a total and multi-valued function defined by 
\begin{equation}\tag{I2}\label{Inp2}{\rm Input}_{\cal M}=\{(\vec x, (1\,.\, \vec \nu\,.\,\bar u))\mid \vec x \in U_{\cal A}^\infty \,\,\&\,\, \vec x \,\genfrac{}{}{0pt}{2}{\longrightarrow} {\,{\rm In}_{\cal M}} (\vec \nu\,.\,\bar u)\} \end{equation} 
where $ \vec x \,\genfrac{}{}{0pt}{2}{\longrightarrow} {\,{\rm In}_{\cal M}} (\vec \nu\,.\,\bar u)$ stands for $(\vec x,(\vec \nu\,.\,\bar u) )\in {\rm In}_{\cal M}$ and $\vec x \,\genfrac{}{}{0pt}{2}{\longrightarrow} {\,{\rm Input }_{\cal M}} (1\,.\,\vec \nu\,.\,\bar u) $ means $(\vec x, (1\,.\, \vec \nu\,.\,\bar u)) \in {\rm Input}_{\cal M}$. Informally, the application of the {\em input procedure} on $\vec x$ can be viewed as a process to provide one of the {\em initial} or {\em start configurations} $con_1\in\{(1\,.\,(n, 1,\ldots, 1)\,.\,\vec x\,.\,\vec y \,.\,(x_n,x_n, \ldots ))\mid \,\vec y\in U_{\cal A}^\infty\}$ by guessing $m\geq 0$ and $m$ individuals $y_1,\ldots,y_m$ which we call {\em guesses}.
A consequence will be that the outputs of a non-deterministic BSS RAM on each input in $U_{\cal A}^\infty$ form a subset belonging to the power set ${\mathfrak P}(U_{\cal A}^\infty)$ of $U_{\cal A}^\infty$. 
\vspace{0.1cm}

\begin{overview}
[Initial values resulting from input and guessing]\label{ConfigurationInfMach}
\hfill

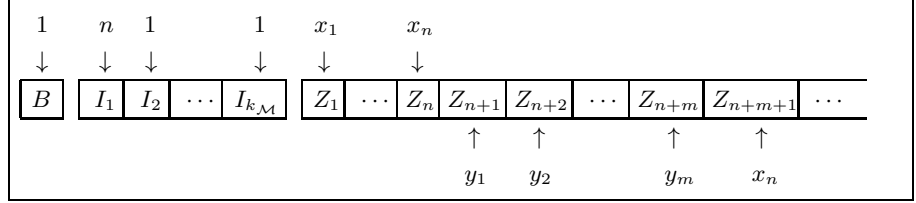
\nopagebreak 
\noindent \fbox{\parbox{11.8cm}{

$\!$\begin{tabular}{l}
{\small
\begin{tabular}{l}
$\!\!1$\\
$\!\!\downarrow$\\ 
\noindent $\!\!\!\! \!\!\! \!$ \begin{tabular}{|c|}
\hline\begin{tabular}{c}$\! \!\!\!\!B\!\!\!\!\! $\\\end{tabular}\\\hline
\end{tabular}\\\\\\
\end{tabular}$\!$
\begin{tabular}{l}
$n$\hspace{0.4cm}$1$\hspace{1.3cm}$1$\\
$\downarrow$\hspace{0.45cm}$\downarrow$\hspace{1.29cm}$\downarrow$\\
\noindent $\!\!\! \!\! $\begin{tabular}{|c|}
\hline
\begin{tabular}{c|}$\!\!\!\! I_1\!\!$\\\end{tabular} 
\begin{tabular}{c|}$\!\!I_2\!\!$\\\end{tabular} 
\begin{tabular}{c|}$\!\!\cdots\!\!$\\\end{tabular} 
\begin{tabular}{c}
$\!\! \!I_{k_{\cal M}} \!\!\! \! \! \!$\\\end{tabular}\\
\hline
\end{tabular}\\\\\\
\end{tabular}$\!$\begin{tabular}{l}

$x_1$\hspace{0.9cm}$x_n$\\
$\,\downarrow$\hspace{1.08cm}$\downarrow$\\
\noindent $\!\!\!$\begin{tabular}{|c}
\hline
\begin{tabular}{c|}$\!\!\!\!\!Z_1\!\!\!$\\\end{tabular}
\begin{tabular}{c|}$\!\!\cdots\!\!$\end{tabular}
\begin{tabular}{c|}$\!\!\!\!Z_n\!\!\!$\\\end{tabular}
\begin{tabular}{c|}$\!\!\!\!Z_{n+1}\!\!\!$\\\end{tabular}
\begin{tabular}{c|}$\!\!\!\!Z_{n+2}\!\!\!$\\\end{tabular}
\begin{tabular}{c|}$\!\!\cdots\!$\end{tabular}
\begin{tabular}{c|}$\!\!\!\!Z_{n+m}\!\!\!$\\\end{tabular}
\begin{tabular}{c|}$\!\!\!\!Z_{n+m+1}\!\!\!\!$\\\end{tabular}
\begin{tabular}{c}$\!\!\cdots\!\!$\end{tabular}\\
\hline
\end{tabular}\\
\hspace{2.05cm}$\uparrow$\hspace{0.7cm}$\uparrow$\hspace{1.63cm}$\uparrow$\hspace{1cm}$\uparrow$\\
\hspace{2cm}$y_1$\hspace{0.55cm}$y_2$\hspace{1.5 cm}$y_m$\hspace{0.75cm}$x_n$\\
\end{tabular}}
\end{tabular}
}}
\end{overview}

\vspace{0.1cm}

The {\em transition system} $({\sf S}_{\cal M}, \to_{\cal M})$ of any non-deterministic BSS RAM ${\cal M}$ that can execute only instructions of the types (1) to (8) is defined as the transition system of any deterministic BSS RAM given in \cite[Definition 4, p.\,584]{GASS20} and, consequently, the binary relations $\to_{\cal M}\,\subset {\sf S}_{\cal M}\times {\sf S}_{\cal M}$ and $(\to_{\cal M})_{{\rm Stop}_{\cal M}}\subset {\sf S}_{\cal M}\times {\sf S}_{\cal M}$ are again single-valued functions on ${\sf S}_{\cal M}$. 

Let ${\sf M}_{\cal A}^{\rm ND}$ be the class of {\sf all these non-deterministic BSS RAMs} over ${\cal A}$ that can only execute instructions of types $(1)$ to $(8)$. 
\begin{definition}[Non-deterministically computable function]\hfill For any\linebreak BSS RAM ${\cal M}$ in ${\sf M}_{\cal A}^{\rm ND}$, its {\em result function} ${\rm Res}_{\cal M}: U_{\cal A}^\infty\to {\mathfrak P}(U_{\cal A}^\infty)$ is defined by \begin{equation}\tag{R2}\label{Res2} {\rm Res_{\cal M}}(\vec x)=\{ {\rm Ouput}_{\cal M}( (\to_{\cal M})_{{\rm Stop}_{\cal M}} (1\,.\,\vec \nu\,.\,\bar u))\mid (\vec x, (\vec \nu\,.\,\bar u))\in {\rm In}_{\cal M}(\vec x)\}
\end{equation} 
for all $\vec x \in U_{\cal A}^\infty$. A partial function $f:\,\subseteq\! U_{\cal A}^\infty\to U_{\cal A}^\infty$ is {\em non-deterministically computable by a BSS RAM over ${\cal A}$} if there is a BSS RAM ${\cal M}$ in ${\sf M}_{\cal A}^{\rm ND}$ whose result function ${\rm Res}_{\cal M}$ assigns to each $\vec x \in U_{\cal A}^\infty$ a set containing at most one tuple that belongs to $U_{\cal A}^\infty$ (which means $|{\rm Res}_{\cal M}({\vec x})|\leq 1$) and the following holds. For $\vec x \in U_{\cal A}^\infty$, ${\rm Res}_{\cal M}({\vec x})=\{f(\vec x)\}$ holds if ${\rm Res}_{\cal M}({\vec x})\not=\emptyset$. $f(\vec x)$ is not defined (denoted shortly by $f(\vec x)\uparrow$) if $\vec x$ is in $U_{\cal A} ^\infty$ and ${\rm Res}_{\cal M}({\vec x})=\emptyset$ holds. A partial multi-valued function $f:\,\subseteq\! U_{\cal A}^\infty \genfrac{}{}{0pt}{2}{\longrightarrow} {\longrightarrow}U_{\cal A}^\infty$ is {\em non-deterministically computable by a BSS RAM over ${\cal A}$} if there is a BSS RAM ${\cal M}$ in ${\sf M}_{\cal A}^{\rm ND}$ such that we have $f=\{(\vec x, \vec z)\mid \vec x \in U_{\cal A}^\infty \,\,\&\,\,\vec z\in {\rm Res}_{\cal M}({\vec x})\}$. $\vec x { \genfrac{}{}{0pt}{2}{\longrightarrow} {f}} \vec z$ means $(\vec x,\vec z)\in f$. 
\end{definition} 

For ${\cal M}\in {\sf M}_{\cal A}^{\rm ND}$ and any input $\vec x\in U_{\cal A}^\infty$, all outputs of ${\cal M}$ are in ${\rm Res}_{\cal M}(\vec x)\subseteq U_{\cal A}^\infty$. We say that ${\cal M}$ {\em halts on $\vec x\in U_{\cal A}^\infty$ for some guesses} if ${\rm Res}_{\cal M}(\vec x)\not=\emptyset$. Then, we also say that ${\cal M}$ {\em accepts} $\vec x$ and we write $ {\cal M} ( \vec x)\downarrow $. Let $H_{\cal M}$ be the {\it halting set of ${\cal M }$} given by $H_{\cal M}=\{ \vec x \in U_{\cal A}^{\infty} \mid {\cal M} ( \vec x)\downarrow \}$ that we also call {\it non-deterministically} {\it semi-decidable} or {\rm ND}-{\it semi-decidable} ({\em by ${\cal M}$}) ({\em over ${\cal A}$}).
Let ${\rm SDEC}_{\cal A}^{\rm ND}$ be the class of {\sf all decision problems ${\rm ND}$-semi-decidable} by a BSS RAM in ${\sf M}_{\cal A}^{\rm ND}$. Thus, $\mathbb{R}_{\geq 0}\in {\rm SDEC}^{\rm ND}_{\mathbb{R}^=}$ holds (cf. Example \ref{PositivRaten}). ${\rm SDEC}_{\cal A}\subseteq {\rm SDEC}_{\cal A}^{\rm ND}$ holds in any case. If $|U_{\cal A}|=1$, then, for any ${\cal M}\in {\sf M}_{\cal A}^{\rm ND}$, ${\rm Input}_{\cal M}$ is a single-valued function and we have ${\rm SDEC}_{\cal A}^{\rm ND}={\rm SDEC}_{\cal A}$. Let us remark that the restriction of the domain of guesses can be useful (cf. Theorem \ref{RestrictGuesses}). If ${\cal A}$ contains two constants $c_1$ and $c_2$, then let $ {\sf M}_{\cal A}^{{\rm DND}}$ be the corresponding class of BSS RAMs ${\cal M}$ with \[{\rm In}_{\cal M}=\{(\vec x,(\,\underbrace{(n, 1,\ldots, 1)}_{\in(\mathbb{N}_+)^{k_{\cal M}}}\,.\,\underbrace{\vec x\,.\,\vec y\,.\,(x_n,x_n, \ldots)}_{\in U_{\cal A}^{\omega}}\,))\mid \vec x\in U_{\cal A}^\infty\,\,\&\,\,\vec y \in \{c_1,c_2\}^\infty\}\] and let ${\rm SDEC}_{\cal A}^{\rm DND}$ be the respective class of {\sf all decision problems} that are {\em digitally non-deterministically semi-decidable} by a machine in ${\sf M}_{\cal A}^{\rm DND}$, and so on.

\subsection{Oracle machines}
The execution of the oracle instructions of type (9) enables the evaluation of queries of the form {\em $(c(Z_1),\ldots, c(Z_{c(I_1)}) \in Q$?} for any {\em oracle} ({\em set}) $Q\subseteq U_{\cal A}^{\infty}$. 
\begin{overview}[Oracle instructions with queries]

\hfill

\nopagebreak 

\noindent \fbox{\parbox{ 11.8cm}{

Oracle instructions

\qquad $(9)$ \quad {\sf $\ell :\,$ if $(Z_1,\ldots,Z_{I_1})\in \!{ \cal O}$ then goto $\ell_1$ else goto $\ell_2$}
}}
\end{overview}
${\cal O}$ stands for any $Q\subseteq U_{\cal A}^\infty$. These oracle instructions are strings and belong to those instructions for which the dots $\ldots$ are part of these strings. It is possible to use oracles such as $O_Q=\bigcup_{n=1}^\infty \{ (z_1,z_2,\ldots, z_{2n} )\in U_{\cal A}^ {2n}\mid (z_2,z_4,\ldots, z_{2n} ) \in Q\}$ and moreover we can use pseudo instructions such as the following instructions.

\vspace{0.1cm}

\qquad$^{\rm pseudo}(9) $\quad {\sf $\ell :\,$ if $(Z_2,Z_4\ldots,Z_{I_1}) \in { \cal O}$ then goto $\ell_1$ else goto $\ell_2$}

\begin{overview}[Further test functions in ${\cal F}_{\cal M}$ for changing the label] \label{Test2}
\hfill

\nopagebreak 
\noindent \fbox{\parbox{11.8cm}{
\begin{itemize}\itemsep0pt

\item 
$T_\ell: \bbbn_+\times U_{\cal A} ^{\omega} \to {\cal L}_ {\cal M}$ \hfill for instructions of type $(9)$
\parskip -0.5mm 
\begin{tabbing} $C_{\ell }(\vec \nu,\bar u)\quad$\=$=( u_1,\ldots,u_{j-1},f_i(u_{j_1},\ldots, u_{j_{m_i}}),u_{j+1}, \ldots)$ \= \kill
$T_{\ell }(\nu_1,\bar u)$\>$=\left\{\begin{array}{ll}\ell_1\quad \mbox{ if $ (u_1,\ldots, u_{\nu_1})\in Q$}\\\ell_2\quad \mbox{ if $ (u_1,\ldots, u_{\nu_1})\not\in Q$}\end{array}\right.$
\end{tabbing} 
\end{itemize}}}
\end{overview}
 The usual {\em oracle BSS RAMs} can only execute instructions of types $(1)$ to $(9)$ and use a set $Q$ as oracle. Let ${\sf M}_{\cal A}^{\big[{\rm ND\big]}}(Q)$ be the class of {\sf all these $\big[$non-$\big]$deterministic oracle BSS RAMs} over ${\cal A}$ using $Q\subseteq U_{\cal A}^{\infty}$ as oracle.\footnote{Note, that, in such descriptions, either only the big square brackets $\big[$ and $\big ]$ or all big square brackets including the enclosed text must be deleted at the same time.} For any ${\cal M} \in {\sf M}_{\cal A}^{\rm ND}(Q)$, ${\rm Input}_{\cal M}$ is defined by (\ref{Inp2}) and, thus, the result function ${\rm Res}_{\cal M}$ assigns a set ${\rm Res}_{\cal M}(\vec x)$ to each input $\vec x$ as defined in (\ref{Res2}). For ${\cal M}\in {\sf M}_{\cal A}(Q)$, ${\rm Input}_{\cal M}(\vec x)$ is defined by (\ref{Inp1}) and ${\rm Res}_{\cal M}(\vec x)$ is defined by (\ref{Res1}). For each ${\cal M}\in {\sf M}_{\cal A}^{\big[{\rm ND\big]}}(Q)$, the {\em transition system} $({\sf S}_{\cal M}, \to_{\cal M})$ is defined analogously to the definition in \cite{GASS20} where $T _{\ell}: (\bbbn\times U_{\cal A}^\omega) \to {\cal L}_{\cal M}$ can, additionally, be used to describe the application of instructions of type (9). Consequently, the relations $\to_{\cal M}$ and $(\to_{\cal M})_{{\rm Stop}_{\cal M}}$ are again functions on ${\sf S}_{\cal M}$. Corresponding to this, we have the classes ${\rm SDEC}_{\cal A}^Q$ and $({\rm SDEC}_{\cal A}^{\rm ND})^Q$ of {\sf all decision problems {\em $\big[{\rm ND}\mbox{-}\big]$semi-decidable}} by a BSS RAM in ${\sf M}_{\cal A}^{\big[{\rm ND\big]}}(Q)$. $P\subseteq U_{\cal A}^\infty$ is {\em decidable over ${\cal A}$ by means of the oracle $Q$} and belongs to ${\rm DEC}_{\cal A}^Q$ if both the set $P$ and its complement $U_{\cal A}^\infty\setminus P$ belong to ${\rm SDEC}_{\cal A}^Q$.

\begin{figure}
\centering
{\begin{tabular}{|l|} \hline\\

\qquad $((\ell \,.\,\vec \nu\,.\,\bar u) , (\,\,\quad\ell_1\,\quad.\,\vec \nu\,.\,\bar u)) \in \,\to_{\cal M}$ \qquad in the case that $(u_1,\ldots, u_{\nu_1}) \in Q$ \\
or \hspace{0.28cm}$((\ell \,.\,\vec \nu\,.\,\bar u), (\,\underbrace{\quad\ell_2\quad}_{ T_{\ell }(\nu_1,\bar u)} .\,\vec \nu\,.\,\bar u))\in \, \to_{\cal M}$ \qquad in the case that $(u_1,\ldots, u_{\nu_1}) \not\in Q$ \\caused by \qquad \qquad \qquad { \sf $\ell \!: \,$ if $(Z_1,Z_2\ldots,Z_{I_1}) \in { \cal O}$ then goto $\ell _1$ else goto 
 $\ell _2$} \\\\\hline \end{tabular}}
\caption{A further pair of configurations $(con_s,con_{s+1})$ in $\to_{\cal M}$ } \label{InstructionsChangeConfInfinite2}
\end{figure}

It is possible to check whether any $x_1\in \mathbb{R}$ is non-negative by an oracle machine in $ {\sf M}_{ \mathbb{R}^=}(\mathbb{R}_{\geq 0})$ (cf. Example \ref{PositivOrakel}). We introduced oracle instruction of type (9) in \cite{GASS97} and used it, for instance, in \cite{GASS07, GASS08A, GASS08B, GASS08C, GASS08D, GASS09A, GASS09B, GASS10}.

\subsection{Other non-deterministic BSS RAMs}

The guessing process by a non-deterministic BSS RAM can be simulated by means of the {\sf non-de\-terministic operator $\nu$} derived from an operator introduced in \cite{MOSCHO}. Let $Q\subseteq U_{\cal A}^\infty$. We will use 
\[\nu[Q](\vec x)=_{\rm df}\{ y_1\in U_{\cal A}\mid ( \vec x \,.\,y_1)\in Q\mbox{ or }(\exists (y_2,\ldots,y_m)\in U_{\cal A}^\infty)(( \vec x \,.\,\vec y) \in Q)\}\] 
for interpreting the following non-deterministic oracle instruction. 

\begin{overview}
[Oracle instructions with Moschovakis' $\nu$-operator]
\hfill

\nopagebreak 

\noindent \fbox{\parbox{ 11.8cm}{

$\nu$-instruction

 \qquad $(10)$ \quad $\ell:\,\, Z_j:=\nu[{\cal O}](Z_1,\ldots,Z_{I_1})$
}}
\end{overview}
In the oracle instructions of type (10), ${\cal O}$ stands for any $Q\subseteq U_{\cal A}^\infty$. The $\nu$-instructions are strings where the dots $\ldots$ are substrings that belong to this version of $\nu$-instructions. The execution of a $\nu$-instruction has the consequence that $ Z_j$ gets any value that belongs to $ {\nu}[Q](c(Z_1),\ldots,c(Z_{n}))$ if $c(I_1)=n$ and $ {\nu}[Q](c(Z_1),\ldots,c(Z_{n}))\not= \emptyset$ hold. Otherwise, this instruction has the effect that the machine loops forever (which means that the stop criterion cannot be reached). 
Any {\em $\nu$-oracle BSS RAM} ${\cal M}$ over ${\cal A}$ is able to execute instructions of types $(1)$ to $(8)$ and (10) and evaluate an oracle $Q\subseteq U_{\cal A}^\infty$ by using $\nu$. ${\rm Input}_{\cal M}$ and ${\rm Output}_{\cal M}$ are again defined by (\ref{Inp1}) and (\ref{Out}). For $Q\subseteq U_{\cal A}^\infty$, let ${\sf M}_{\cal A}^\nu(Q)$ be the class of {\sf all $\nu$-oracle BSS RAMs} over ${\cal A}$ that can evaluate $Q$ by using $\nu$.

\begin{overview}[Multi-valued maps in ${\cal F}_{\cal M}$ for changing a state]\label{Oper2}
\hfill

\nopagebreak 
\noindent \fbox{\parbox{11.8cm}{
\begin{itemize}\itemsep0pt
\item
$F_\ell:\,\subseteq\!(\bbbn_+\times U_{\cal A} ^{\omega})\genfrac{}{}{0pt}{2}{\longrightarrow} {\longrightarrow} U_{\cal A} ^{\omega}$ \hfill for instructions of type $(10)$
\parskip -0.5mm 
\begin{tabbing} 
$C_{\ell }(\vec \nu,\bar u)\quad$\=$=( u_1,\ldots,u_{j-1},f_i(u_{j_1},\ldots, u_{j_{m_i}}),u_{j+1}, \ldots)$ \= \kill
$F_\ell=\{((\nu_1,\bar u),(u_1,\ldots,u_{j-1},y , u_{j+1}, \ldots))\,\,\mid\,\, y\in \nu[Q] (u_1,\ldots, u_{\nu_1})\}$
\end{tabbing} 
\end{itemize}}}
\end{overview}

For any BSS RAM ${\cal M}$ in ${\sf M}_{\cal A}^\nu(Q)$, the relations $\to_{\cal M}$ and $(\to_{\cal M})_{{\rm Stop}_{\cal M}}$ are {\sf partial multi-valued functions} $\to_{\cal M}:\,\,\subseteq\!{\sf S}_{\cal M} \genfrac{}{}{0pt}{2}{\longrightarrow} {\longrightarrow} {\sf S}_{\cal M}$ and $(\to_{\cal M})_{{\rm Stop}_{\cal M}}:\,\,\subseteq\!{\sf S}_{\cal M} \genfrac{}{}{0pt}{2}{\longrightarrow} {\longrightarrow} {\sf S}_{\cal M}$. Each configuration of ${\cal M}$ assigned to a configuration $(\ell\,.\,\vec \nu\,.\,\bar u)$ by $\to_{\cal M}$ is uniquely defined {\sf or} in $\{(\ell+1\,.\,\vec \nu\,.\, (u_1,\ldots, u_{j-1},y,u_{j+1},\ldots) )\mid y\in \nu[Q] (u_1,\ldots, u_{\nu_1})\}$. A consequence will be that we also write $\genfrac{}{}{0pt}{2}{\longrightarrow} {\longrightarrow}_{\cal M}$ instead of $\to _{\cal M}$. According to this, the multi-valued function $(\to_{\cal M})_{{\rm Stop}_{\cal M}}$, that we also denote by $(\genfrac{}{}{0pt}{2}{\longrightarrow} {\longrightarrow}_{\cal M})_{{\rm Stop}_{\cal M}}$, is the {\sf intersection} of the transitive closure of $\genfrac{}{}{0pt}{2}{\longrightarrow} {\longrightarrow}_{\cal M}$ and the set ${\cal S}_{\cal M}\times \{(\ell_{{\cal P}_{\cal M}}\,.\,\vec \nu \,.\, \bar u)\mid \vec \nu \in (\bbbn_+)^{k_{\cal M}}\,\,\&\,\, \bar u\in U_{\cal A}^\omega\} $. Consequently, by (\ref{Res3}), ${\rm Res}_{\cal M}(\vec x)$ is a set and all outputs of ${\cal M}$ for the input $\vec x$ belong to ${\rm Res}_{\cal M}(\vec x)$. 

\begin{definition}[$\nu$-computable functions] Let $Q\subseteq U_{\cal A}^\infty$ be given. For every $\nu$-oracle BSS-RAM ${\cal M}$ in ${\sf M}_{\cal A}^\nu(Q)$, the {\em the result function} ${\rm Res}_{\cal M}: U_{\cal A}^\infty\to {\mathfrak P}(U_{\cal A}^\infty)$ is defined by 
\begin{equation}\tag{R3}\label{Res3} {\rm {\rm Res}_{\cal M}}({\vec x})=\{{\rm Output}_{\cal M}(con)\mid ({\rm Input}_{\cal M}(\vec x),con)\in (\genfrac{}{}{0pt}{2}{\longrightarrow} {\longrightarrow}_{\cal M})_{{\rm Stop}_{\cal M}}\}\end{equation} 
for all $\vec x \in U_{\cal A}^\infty$. A partial function $f:\,\,\subseteq\! U_{\cal A}^\infty\to U_{\cal A}^\infty$ is non-deterministically {\em $\nu$-computable by a BSS RAM over ${\cal A}$} ({\em by means of the oracle $Q$}) if there is a BSS RAM ${\cal M}$ in ${\sf M}_{\cal A}^\nu(Q)$ such that ${\rm Res}_{\cal M}({\vec x})=\{f(\vec x)\}$ holds if ${\rm Res}_{\cal M}({\vec x})\not=\emptyset$ and $f(\vec x)\uparrow$ holds if $\vec x \in U_{\cal A}^\infty$ and ${\rm Res}_{\cal M}({\vec x})=\emptyset$ hold. A partial multi-valued function $f:\,\,\subseteq\! U_{\cal A}^\infty \genfrac{}{}{0pt}{2}{\longrightarrow} {\longrightarrow}U_{\cal A}^\infty$ is non-deterministically {\em $\nu$-computable by a BSS RAM over ${\cal A}$} ({\em by means of $Q$}) if there is a BSS-RAM ${\cal M}$ in ${\sf M}_{\cal A}^\nu(Q)$ such that $f=\{(\vec x, \vec y)\mid \vec x \in U_{\cal A}^\infty \,\,\&\,\,\vec y\in {\rm Res}_{\cal M}({\vec x})\}$ holds. Let $ \vec x \genfrac{}{}{0pt}{2}{\longrightarrow} {f} \vec y$ again mean $(\vec x,\vec y)\in f$. 
\end{definition} 
Each halting set $H_{\cal M}$ of a machine ${\cal M}\in {\sf M}_{\cal A}^\nu(Q)$ is called ({\it non-deterministically}) {\it $\nu$-semi-decidable by ${\cal M}$ over ${\cal A}$}.
Let $({\rm SDEC}_{\cal A}^\nu)^Q$ be the class of {\sf all problems} that are {\sf $\nu$-semi-decidable} by a machine in ${\sf M}_{\cal A}^\nu(Q)$. The ideas leading to Proposition 10 in \cite[p.\,592]{GASS20} (and \cite[Proposition 1]{GASS19Preprint}) can be used to prove the following statements which include a result presented in 2015 (for details see \cite{GaVV15}). For the proofs, see Part II.
\begin{theorem}\label{OracleMNondetM}Let ${\cal A}$ be a first-order structure of any signature. If $Q\in {\rm SDEC}_{\cal A}$, then ${\rm SDEC}_{\cal A}\subseteq({\rm SDEC}_{\cal A}^\nu)^Q\subseteq {\rm SDEC}_{\cal A}^{\rm ND}$. 
 \end{theorem}
\begin{theorem}Let ${\cal A}$ be a first-order structure of any signature. If $Q=U_{\cal A}^\infty$, then $({\rm SDEC}_{\cal A}^\nu)^Q={\rm SDEC}_{\cal A}^{\rm ND}$. 
\end{theorem}
It is also possible to restrict the considered oracle set to sets of constants. 
\begin{theorem}\label{RestrictGuesses}Let ${\cal A}$ be any first-order structure with two constants $c_1$ and $c_2$ and let $Q=\{c_1,c_2\}^\infty$. If $Q\in {\rm SDEC}_{\cal A}$, then 
$({\rm SDEC}_{\cal A}^\nu)^Q={\rm SDEC}_{\cal A}^{\rm DND}\subseteq{\rm SDEC}_{\cal A}^{\rm ND}$.
\end{theorem}
Let us remark that this and other deterministic and non-deterministic operators were discussed at several meetings and the discussions are summarized in several proceedings. Some of the considerations regarding operators such as Moschovakis' operator can be found in \cite{GaVV15,GASS16,GaPaSt17,GaPaSt18,GaPaSt21}.

\vspace{0.1cm}

If no constants are available, then we could use the following type of instructions for guessing labels by {\em non-deterministic branching}.

\begin{overview}
[Instructions for non-deterministic branching]

\hfill

\nopagebreak 

\noindent \fbox{\parbox{ 11.8cm}{

Instructions for guessing labels

$(11)$ \quad {\sf $\ell :\,$ goto $\ell_1$ or goto $\ell_2$}
}}
\end{overview}
This type of instructions is known from discussions by Armin Hemmerling \cite{{H98}} and others. The non-deterministic branching instructions are also used in \cite{GASS01}.

\begin{overview}[A multi-valued map in ${\cal F}_{\cal M}$ for changing the label]
\hfill

\nopagebreak 
\noindent \fbox{\parbox{11.8cm}{
\begin{itemize}\itemsep0pt
\item $G: {\cal L}_{\cal M} \genfrac{}{}{0pt}{2}{\longrightarrow} {\longrightarrow} {\cal L}_{\cal M}$ \hfill for instructions of type $(11)$
\parskip -0.5mm 
\begin{tabbing} 
$C_{\ell }(\vec \nu,\bar u)\quad$\=$=( u_1,\ldots,u_{j-1},f_i(u_{j_1},\ldots, u_{j_{m_i}}),u_{j+1}, \ldots)$ \= \kill
\quad $\ell \genfrac{}{}{0pt}{2}{\longrightarrow} {G} \ell_0$ \hspace*{6.7cm} with $\ell_0\in \{\ell_1, \ell_2\}$\\which means\\
\quad $G=\{(\ell,\ell_1), (\ell,\ell_2)\}$ and $G_\ell=\{\ell_1,\ell_2\}$

\end{tabbing} 
\end{itemize}}}
\end{overview}

\begin{figure}
\centering
{\begin{tabular}{|l|} \hline\\
\qquad \quad $((\ell \,.\,\vec \nu\,.\,\bar u) , (\hspace{0.23cm}\ell_1\hspace{0.23cm} .\,\vec \nu\,.\,\bar u)) \in \,\to_{\cal M}$ \hspace{4.2cm} \hfill if $\ell \genfrac{}{}{0pt}{2}{\longrightarrow} {G} \ell_1$ \\
and \quad $((\ell \,.\,\vec \nu\,.\,\bar u), (\,\underbrace{\ell_2}_{ \in G_{\ell }} .\,\vec \nu\,.\,\bar u))\in \, \to_{\cal M}$ \hfill if $\ell \genfrac{}{}{0pt}{2}{\longrightarrow} {G} \ell_2$ \\caused by \hfill { \sf $\ell :\,$ goto $\ell_1$ or goto $\ell_2$} \\\\\hline 
\end{tabular}}
\caption{Two pairs of configurations $(con_s,con_{s+1})$ in $\to_{\cal M}$ } \label{InstructionsChangeConfInfinite3}
\end{figure}

Let ${\sf M}_{\cal A}^{\rm NDB}$ be the class of {\sf all non-deterministic BSS RAMs} over ${\cal A}$ that can only execute instructions of types $(1)$ to $(8)$ and (11). For any BSS RAM ${\cal M}$ in ${\sf M}_{\cal A}^{\rm NDB}$ the corresponding relations $\to_{\cal M}$ and $(\to_{\cal M})_{{\rm Stop}_{\cal M}}$ are multi-valued functions denoted also by $\genfrac{}{}{0pt}{2}{\longrightarrow} {\longrightarrow}_{\cal M}$ and $(\genfrac{}{}{0pt}{2}{\longrightarrow} {\longrightarrow}_{\cal M})_{{\rm Stop}_{\cal M}}$. $\genfrac{}{}{0pt}{2}{\longrightarrow} {\longrightarrow}_{\cal M}$ defined now for the instructions of types (1) to (8) and (11) is a total multi-valued function and $(\genfrac{}{}{0pt}{2}{\longrightarrow} {\longrightarrow}_{\cal M})_{{\rm Stop}_{\cal M}}$ is a partial multi-valued function. ${\rm Input}_{\cal M}$ and ${\rm Output}_{\cal M}$ can be defined by using (\ref{Inp1}) and (\ref{Out}). Consequently, ${\rm Res}_{\cal M}(\vec x)$ should be again a set and the outputs of such a machine ${\cal M}$ in ${\sf M}_{\cal A}^{\rm NDB}$ on $\vec x$ should belong to ${\rm Res}_{\cal M}(\vec x)$. 

\begin{definition}
[Computability by non-deterministic branching]\hfill For \linebreak every BSS-RAM ${\cal M}$ in ${\sf M}_{\cal A}^{\rm NDB}$, the function ${\rm Res}_{\cal M}: U_{\cal A}^\infty\to {\mathfrak P}(U_{\cal A}^\infty)$ that can be defined for all $\vec x \in U_{\cal A}^\infty$ by {\rm (\ref{Res3})} is called {\em the result function of ${\cal M}$}. 
\end{definition} 
Let ${\rm SDEC}_{\cal A}^{\rm NDB}$ be the class of {\sf all decision problems non-deterministically branching semi-decidable} by a BSS RAM in ${\sf M}_{\cal A}^{\rm NDB}$.

\begin{theorem}Let ${\cal A}$ be any first-order structure containing the identity relation and two constants $c_1$ and $c_2$ and let $Q=\{c_1,c_2\}^\infty$. If $Q\in {\rm SDEC}_{\cal A}$, then ${\rm SDEC}_{\cal A}\subseteq ({\rm SDEC}_{\cal A}^\nu)^Q={\rm SDEC}_{\cal A}^{\rm DND}={\rm SDEC}_{\cal A}^{\rm NDB}$.
\end{theorem}
\section{Examples}\label{SectionExamples}

\subsection{Input and output functions}

Some examples for infinite ${\cal A}$-machines that are not BSS RAMs are given in \cite{GASS20}. For instance, it is possible that ${\cal A}_0$-machines compute functions of the form $f:\,\subseteq\!\bbbn_+\to \bbbn_+$ and $f:\,\subseteq\!\{0,1\}^* \to \{0,1\}^*$, respectively, by simulating Turing machines. Here, $\{0,1\}^*$ is the set of all strings (or words) over $\{0,1\}$. 

\begin{example}[Input and output functions for machines over ${\cal A}_0$]\label{Example0TM}\hfill 
Let \linebreak ${\rm bin}:\mathbb{N}_+\to \{0,1\}^*$ be defined by ${\rm bin}(\sum_{i=1}^n x_i\cdot 2^{i-1})=x_n\cdots x_1$ for all $\vec x \in\{0,1\}^\infty $ with $x_n=1$ and let ${\rm bin}^{-1}$ be the partial inverse of ${\rm bin}$.
Machines of the form $(\{0,1\}^\omega,(\mathbb{N}_+)^k, {\cal L}, {\cal P}, {\cal A}_0, {\rm In}, {\rm Out})$ can contain different input and output functions ${\rm In}_*$ and ${\rm Out}_*$ or ${\rm In}_{\mathbb{N}}$ and ${\rm Out}_{\mathbb{N}}$ defined as follows, and so on (cf.\,\cite{GASS20}).

\begin{center}
\fbox{\begin{minipage}{11.7cm} 

\vspace{0.1cm}

${\rm In}_{\mathbb{N}}(m)={\rm In}_*\circ{\rm bin}(m)$ \hfill for $m\in\mathbb{N}_+$

${\rm In}_*(x_n\cdots x_1)=((2n,1,\ldots,1)\,.\,(0,x_1,0,x_2,\ldots,0,x_n,1,1,1,1,\ldots))$

 \hfill with $(2n,1,\ldots,1)\in (\mathbb{N}_+)^k$

\vspace{0.1cm}
 
${\rm Out}_{\mathbb{N}} (\nu_1,\ldots,\nu_k,u_1,u_2,u_3,\ldots)= {\rm bin}^{-1}\circ{\rm Out}_* (\nu_1,\ldots,\nu_k,u_1,u_2,u_3,\ldots)$ 

${\rm Out}_* (\nu_1,\ldots,\nu_k,u_1,u_2,u_3,\ldots)=\Lambda$ \hfill if $u_{2\lambda_0-1}=1$ (for $\lambda_0$ and $\lambda_1$, see \,\cite{GASS20})

${\rm Out}_* (\nu_1,\ldots,\nu_k,u_1,u_2,u_3,\ldots)=u_{2\lambda_0}u_{2\lambda_0+2}\cdots u_{2\lambda_1}$ \hfill otherwise
\end{minipage}}
\end{center}
\end{example}

\subsection{Computing functions and deciding sets}
\begin{example}[Computing the sum of three numbers]\label{Example1} The computation of $\sum_{i=1}^3 x_i$ can be done by a finite $\mathbb{R}^=$-machine ${\cal M}$ using three registers $Z_1$, $Z_2$, and $Z_3$ and the following program that belongs to ${\sf P}_{(2;2,2,2;2)}$. 

\vspace{0.1cm}

{\sf

$1: \, Z_1:=f_1^2 (Z_1,Z_2);$

$2: \, Z_1:=f_1^2 (Z_1,Z_3);$

$3:\,$ {\sf stop.}}

\noindent The input function is given by ${\rm In}_{\cal M}( x_1,x_2,x_3)=(x_1,x_2,x_3)$ and the output function is defined by ${\rm Out}_{\cal M}(u_1,u_2,u_3)=u_1$.

\end{example}
As mentioned above, for checking equations such as $\sum_{i=1}^{100} x_i=1$ and $\sum_{i=1}^{1000} x_i=1$, it would be better to use another type of machines.
\begin{example}[Deciding whether the sum of $n$ numbers is 1]\label{VieleSummanden}\hfill Whether \linebreak the sum $\sum_{i=1}^{n} x_i$ of complex numbers in $(\mathbb{C};1,0,i;+,-,\cdot;=)$ is equal to $1$, can be checked by a BSS RAM ${\cal M}$ by using the following $(3;2,2,2;2)$-program. Let ${\cal M}\in {\sf M}_{\mathbb{C}}$ be equipped with three index registers $I_1$, $I_2$, and $I_3$. Thus, for any input $\vec x\in \mathbb{C}^n$, the start configuration is $(1\,.\,(n,1,1)\,.\,(x_1,\ldots, x_n,x_n,x_n,\ldots))$. We have, for example, ${\rm Input}_{\cal M}(3i,2,3,6+i, 7)= (1\,.\,(5,1,1)\,.\,(3i,2,3,6+i,7,7,7,\ldots)) $. 

\vspace{0.1cm}

{\sf 

$1: \, I_2:=I_2+1;$

$2 :$ {\sf if $I_1=I_3$ then goto $7$ else goto $3$}$;$ 

$3: \, I_3:=I_3+1;$

$4: \, Z_{I_2}:=Z_{I_3};$

$5: \, Z_1:=f_1^2 (Z_1,Z_2);$

$6 :$ {\sf if $I_1=I_3$ then goto $7$ else goto $3$}$;$ 

$7$: $Z_2:=c_1^0;$

{\sf $8: \,$ if $r_1^{2}(Z_1,Z_2)$ then goto $10$ else goto $9$}$;$

$9:\, Z_1:=c_2^0;$

$10:\, I_1:=1;$

$11:$\, {\sf stop.}

\vspace{0.1cm}
}
\noindent At the beginning, $I_1$ contains $n$ and the content of $I_2$ and $I_3$ is 1. $I_3$ is used for storing the index $i$ of the currently considered summand $x_i$. The instructions 1 to 6 allow to compute the sum. If ${\cal M}$ halts on an input and the output is 1, then the input is accepted. Otherwise, the output is 0 and the input is rejected. 
\end{example}

\begin{example}[${\rm ND}$-semi-deciding $\mathbb{R}_{\geq 0}$]\label{PositivRaten} The set $\mathbb{R}_{\geq 0}$ of all non-negative reals can be described by $\mathbb{R}_{\geq 0}=\{x_1\in\mathbb{R} \mid (\exists y_1\in \mathbb{R})(x_1=y_1 \cdot y_1)\}$. To check whether $x_1$ is non-negative by a machine ${\cal M}\in {\sf M}^{\rm ND}_{\mathbb{R}^=}$, the following $(2;2,2,2;2)$-program ${\cal P}$ can be used if ${\cal M}$ is equipped with two index registers $I_1$ and $I_2$. If $x_1$ is the input, then we obtain, for any $m\geq 1$ and any guesses $y_1,\ldots,y_m$, a start configuration $(1\,.\,(1,1)\,.\,(x_1,y_1,\ldots,y_m,x_1,x_1,\ldots))\genfrac{}{}{0pt}{2}{\longleftarrow} {\, {\rm Input}_{\cal M}} x_1$ defined by (\ref{Inp2}). 

\vspace{0.05cm}

 {\sf $1: \,$ if $I_1=I_2$ then goto $2$ else goto $1$}$;$

 {\sf $2: \, Z_3:=f_3^2 (Z_2,Z_2)$}$;$

 {\sf $3: \,$ if $r_1^{2}(Z_1,Z_3)$ then goto $4$ else goto $3$}$;$

 {\sf $4:\,$ stop.}

\vspace{0.05cm}

 \noindent At the beginning, an input $x_1\in \mathbb{R}$ of length 1 is stored in $Z_1$ and the content of $I_1$ is then 1. In this case, the first guessed number is stored in $Z_2$. ${\cal M}$ semi-decides $\mathbb{R}_{\geq 0}$ non-deterministically by executing ${\cal P}$.
\end{example}

\begin{example}[Semi-deciding $\mathbb{R}_{\geq 0}$ and $\mathbb{R}^\infty \setminus\mathbb{R}_{\geq 0}$ by an oracle]\label{PositivOrakel} To semi-decide $\mathbb{R}_{\geq 0}$, the following program and the oracle set $\mathbb{R}_{\geq 0}$ can be used. 

\vspace{0.05cm}

 {\sf $1: \,$ if $I_1=I_2$ then goto $2$ else goto $1$}$;$

 {\sf $2: \,$ if $(Z_1,\ldots,Z_{I_1})\in \!{ \cal O}$ then goto $3$ else goto $2$}$;$

 {\sf $3:\,$ stop.}

\vspace{0.05cm}

\noindent The complement $\mathbb{R}^\infty \setminus\mathbb{R}_{\geq 0}$ is semi-decidable by using the following program. 

\vspace{0.05cm}

 {\sf $1: \,$ if $I_1=I_2$ then goto $2$ else goto $3$}$;$

 {\sf $2: \,$ if $(Z_1,\ldots,Z_{I_1})\in \!{ \cal O}$ then goto $2$ else goto $3$}$;$

 {\sf $3:\,$ stop.}

\vspace{0.05cm}

\noindent Thus, $\mathbb{R}_{\geq 0}$ can also be decided by a single oracle BSS RAM ${\cal M}$ in $ {\sf M}_{\mathbb{R}^=}(\mathbb{R}_{\geq 0})$:

\vspace{0.05cm}

{\sf $1: \,$ $Z_2:=c_2^0;$

$2: \,$ $I_2=I_2+1;$

$3: \,$ if $(Z_1,\ldots,Z_{I_1})\in \!{ \cal O}$ then goto $4$ else goto $5;$

$4:\,$ $Z_2:=c_1^0;$

$5:\, $ $ I_1:=1;$

$6:\, $ $Z_{I_1}:= Z_{I_2};$

$7:$\, stop.}

\vspace{0.05cm}

\noindent Since $\,{\rm Input}_{\cal M}(4,6, 7)= (1\,.\,(3,1)\,.\,(4,6,7,7,\ldots) )$ $\to _{\cal M}(\textcolor{red}{2},3,1,4,\textcolor{blue}{0},7,7,\ldots)\to _{\cal M}(\textcolor{red}{3},3,\textcolor{blue}{2},4,0,7,7,\ldots)\to _{\cal M}(\textcolor{red}{5},3,2,4,0,6,7,7,\ldots)$ $ \to _{\cal M}(\textcolor{red}{6},\textcolor{blue}{1},2,4,0,6,7,7,\ldots)$\linebreak $\to _{\cal M}(\textcolor{red}{7},1,2,\textcolor{blue}{0},0,6,7,7,\ldots)=(\to_{\cal M})_{{\rm Stop}_{\cal M}}({\rm Input}_{\cal M}(4,6, 7))$ holds, we obtain, for example, ${\rm Res}_{\cal M}(4,6, 7)={\rm Output}_{\cal M}(7,1,2,0,\ldots) ={\rm Out}_{\cal M}(1,2,0,\ldots)=0$. 
\end{example}

\subsection{Evaluating formulas}
\begin{example}[Checking a formula by a finite ${\cal A}$-machine]\label{Example2}$\!$Let $U$ be a set of sets. In the context of computation, we use the identity relation $=$ like any other relation over $U$ if it should be possible to decide the equality of any two sets in $U$ by executing one single instruction. In line with Agreement \ref{AgreeIdent}, let ${\cal A}$ be the structure $(U;\emptyset;; \in,=)$ of signature $\sigma=_{\rm df}(1;; 2,2)$. In this case, the relation $=$ allows to interpret the symbol $r_2^2$ whenever $r_2^2$ occurs in a formula or a program of signature $\sigma$. Let $\phi(X_1, \ldots,X_n)$ be a quantifier-free first-order formula of signature $\sigma$ that can contain the symbol $c^0_1$ for a constant (here for $\emptyset$), the symbols $r_1^2$ and $r_2^2$ (or, for simplicity matters, we could use $=$ instead of $r_2^2$) for binary relations such as $\in$ and $=$, and the variables $X_1, \ldots,X_n$. Let ${\cal P}_\phi$ be a $\sigma$-program for evaluating the formula $\phi$ by a finite ${\cal A}$-machine that is equipped with a fixed number $j_{\cal M}$ of registers $Z_1,\ldots,Z_{n+1}$ ($j_{\cal M}=n+1$) for storing the input $(x_1,\ldots,x_n)$ and the constant $\emptyset$ such that ${\cal P}_\phi$ halts on $(x_1, \ldots,x_n)$ if and only if $\phi(X_1, \ldots,X_n)$ becomes true when all $x_i$ ($i\in \{1,\ldots,n\}$) are assigned to $X_i$ and the symbols $c_1^0$, $r_1^2$, and $r_2^2$ are interpreted by $\emptyset$, $\in$, and $=$. Here, we can use the fact that every propositional formula (that can be derived from $\phi$ by replacing every atomic formula by a new propositional variable) is equivalent to a formula that is in disjunctive normal form. Therefore, for any input, it is sufficient to evaluate $\phi$ in the same way as a formula in disjunctive normal form by using the truth values of its atomic subformulas of the form $r_i^2(X_{j_1},X_{j_2})$, $r_i^2(X_{j_1},c_1^0)$, $\ldots$ in $\phi(X_1, \ldots,X_n)$. These formulas can be evaluated by branching instructions whose conditions have the form $r_i^2(Z_{j_1},Z_{j_2})$ and imply which instruction follows depending on the truth-values obtained. The stop-instruction should be reached only if one of the conjunctions in the corresponding normal form is true. Thus, $\{(x_1, \ldots,x_n)\mid \phi(x_1, \ldots,x_n) \mbox{ holds in } {\cal A}\}$ is semi-decidable over ${\cal A}$. 
\end{example}

Note, that the evaluation of a formula by a finite machine can, in individual cases, be faster performed when the polish notation of such a formula is used. 

If we want to check all first-order quantifier-free formulas of any signature $\sigma$ (for finite sets $N_2$ and $N_3$) by one machine, we need suitable codes of the formulas. We know that these formulas, their polish notations, and their truth values can be defined recursively (see e.g.\,\,also \cite[pp.\,\,22--26]{A75}). Consequently, the evaluation of such formulas can be done by using $Z$-registers like a stack. However, the evaluation of first-order formulas of signature $\sigma$ can be only done by a machine over a structure ${\cal A}$ that allows to encode and to evaluate the formulas. Moreover, it can be necessary to use an infinite machine, if these formulas or the finite sequences of individuals stored in the stack cannot be encoded by a finite number $j$ of individuals for a fixed $j$. If ${\cal A}$ contains only a finite number of individuals, a finite ${\cal A}$-machine is, in general, not sufficient. Example \ref{EvalMitStack} demonstrates the use of a stack given by $Z$-registers of an infinite ${\cal A}$-machine ${\cal M}$. Its halting set (that is {\em semi-decidable over ${\cal A}$}) is a subset of ${\sf I}_{\cal M}$. 
\begin{example}[Evaluation of Boolean formulas by an ${\cal A}$-machine]\label{EvalMitStack}\hfill
Let \linebreak $\phi$ be a first-order Boolean formula of signature $(2;1,2,2,2;)$ without variables that can be interpreted by using $(\{0,1\};0,1;\neg,\land,\lor,\leftrightarrow;)$ defined completely as follows. In the following, we do not differentiate between the constants and operations themselves, on the one hand, and their symbols, on the other hand.
\begin{center}
\begin{tabular}{c||c|c|c}
$x\quad \quad y$&$x\land y$&$x\lor y$&$\leftrightarrow$\\
\hline
0\quad\quad 0&0&0&1\\
0\quad\quad 1&0&1&0\\
1\quad\quad 0&0&1&0\\
1\quad\quad 1&1&1&1\\
\end{tabular} \hspace{2cm}
\begin{tabular}{c||c}
$x$&$\neg x$\\
\hline
0&1\\
1&0\\
\end{tabular}
\end{center}

\noindent
We obtain the same truth values for any formula in infix and prefix notation regardless of the notation used. $0$ and $1$ are symbols and the truth values themselves, respectively. The corresponding holds for $\neg,\land,\lor,\leftrightarrow$. $=$ is here a relation between truth values allowing to say (indirectly) that two formulas are equivalent.

\vspace{0.2cm}

\noindent \begin{tabular}{ll}
$\underbrace{(\textcolor{red}{(0\land 1)} \lor \neg\textcolor{red}{(1\leftrightarrow 0)} )\land ( \neg \textcolor{red}{(0\land 1)} \lor \textcolor{red}{(0 \leftrightarrow 0)} )} _ {\mbox{\small infix notation}}$ 
& $\underbrace{\land \lor \land 01 \neg\leftrightarrow 10 \lor \neg \land 01 \textcolor{red}{\leftrightarrow 00} } _ {\mbox{\small prefix notation}}$\\\vspace{0.2cm}
$= \textcolor{red}{(0 \lor \neg 0 )}\land \textcolor{red}{( \neg 0 \lor 1 )}$ 
& $= \land \lor \land 01 \neg\leftrightarrow 10 \lor \neg \,\textcolor{red}{\land 01} 1$\\
$=\textcolor{red}{1}\land \textcolor{red}{1} $
&$= \land \lor \land 01 \neg\leftrightarrow 10 \lor \textcolor{red}{\neg 0} 1$\\
$=1$
&$= \land \lor \land 01 \neg\leftrightarrow 10 \, \textcolor{red}{\lor\, 11}$\\
&$= \land \lor \land 01 \neg\, \textcolor{red}{\leftrightarrow 10} 1$\\
&$= \land \lor \land 01 \textcolor{red}{\neg 0} 1$\\
&$= \land \lor \textcolor{red}{\land 01} 1 1$\\
&$= \land \,\textcolor{red}{\lor \,0 1} 1$\\
&$= \textcolor{red}{\land 1 1}$\\
&$= 1$\\
\end{tabular}

\vspace{0.1cm}

\noindent
For simplicity matters, let ${\cal M}$ be an infinite $(U; U;;=)$-machine with $U=\{0,1,\neg,\land,\lor,\leftrightarrow,\land_0,\lor_0,\leftrightarrow_0, \land_1,\lor_1,\leftrightarrow_1,\Lambda,\# ,q_0\}$ to evaluate Boolean formulas without variables. Each individual is a constant that can be used in executing instructions of type (2) and the identity can be used in branching instructions. The input space ${\sf I}_{\cal M}$ is the set of strings given by ${\sf I}_{\cal M}=\{1,0,\neg,\land,\lor,\leftrightarrow,(,)\}^*$. If $s_1\cdots s_n\in {\sf I}_{\cal M}$ is a first-order Boolean formula $\phi$ without variables (given in the usual infix notation), then let $(u_1, \ldots,u_m)$ be the tuple $\,^{\rm pol}\phi\in\! \{0,1,\neg,\land,\lor,\leftrightarrow\}^\infty$\linebreak that corresponds to the polish prefix notation of $\phi$ in $ \{0,1,\neg,\land,\lor,\leftrightarrow\}^*$ and let \[{\rm In}_{\cal M}(s_1\cdots s_n)=((m,1,\ldots,1)\,.\,(u_1, \ldots,u_m,u_m,\ldots)).\] Otherwise, let ${\rm In}_{\cal M}(s_1\cdots s_n)=((1,1,\ldots,1)\,.\,( 0, 0, 0,\ldots))$. In this way, we get, for example, ${\rm In}_{\cal M}((1\land0)\land1\lor1)=((1,1,\ldots,1)\,.\, ( 0, 0, 0,\ldots)))$ and
 \[\begin{array}{ll}{\rm In}_{\cal M}(\,( \overbrace{\,\overbrace{(\textcolor{blue}{1}\textcolor{red}{\land}\textcolor{magenta}{0}\,)}^ {\phi_1}\textcolor{brown}{\lor}\textcolor{green}{1})}^ {\phi}\,)&=((\overbrace{5,1,\ldots,1}^{\rm indices})\,.\, \overbrace{^{\rm pol}\phi \,.\,(\underbrace{1,1,1,1,1,1,\ldots\,\qquad}_{\rm any \, 1 \, is \, the\, last\, input\, value\,\textcolor{green}{1}})}^{\rm individuals })\vspace{0.1cm}\\
&=((5,1,\ldots,1)\,.\, \underbrace{(\textcolor{brown}{\lor}\,.\, ^{\rm pol}\phi_1\,.\,\textcolor{green}{1})}_{ ^{\rm pol}\phi }\,.\,(1,1,\ldots))\\
\vspace{0.05cm}
&=((5,1,\ldots,1)\,.\, (\underbrace{\textcolor{brown}{\lor},\underbrace{\textcolor{red}{\land}, \textcolor{blue}{1},\textcolor{magenta}{0}}_{ ^{\rm pol}\phi_1},\textcolor{green}{1}}_{ ^{\rm pol}\phi },1,1,\ldots)).\end{array}\]
Let ${\rm Input}_{\cal M}({\sf i})= (1\,.\,{\rm In}_{\cal M}({\sf i}))$ for all ${\sf i}\in {\sf I}_{\cal M}$.

To evaluate a formula given by $ ^{\rm pol}\phi $, ${\cal M}$ can work like a pushdown automaton $N$ (cf.\,\cite{Schoening08}). The individuals $\land_0,\lor_0,\leftrightarrow_0, \land_1,\lor_1,\leftrightarrow_1$ are used like other states of $N$. They allow to say which of the unary operations defined as follows is to be executed.

\begin{center}
\begin{tabular}{c||c|c|c|c|c|c}
$x$&$\land_1 x$&$\land_0x$&$\lor_1x$&$\lor_0x$&$\leftrightarrow_1x$&$\leftrightarrow_0x$\\
\hline
0&0&0&1&0&0&1\\
1&1&0&1&1&1&0\\
\end{tabular}\end{center}

\noindent At the beginning, $q_0$ corresponds to the initial state of $N$. $q_0$ says that the last symbol of the rest of the input in polish notation has to be evaluated. $\Lambda$ says that the evaluation of the formula is finished. The evaluation of the formula can be performed by means of a stack which is represented by a tuple $(c(Z_i),\ldots, c(Z_j))$ with $i\geq j$ such that $c(Z_i)$ is the last symbol stored in the stack and $Z_j$ generally contains $\#$. $\#$ is called the lowest stack symbol. If $i>j$, then we have $c(Z_j)=\#$ and the other components $c(Z_i),\ldots, c(Z_{j+1})$ are not equal to $\#$. The last stored symbol $c(Z_i)$ is the first individual that can be read and removed from the stack. If $c(Z_i) =\#$, then $i=j$ holds which means that no important information is in the stack. Otherwise, without the symbol $\#$ in the stack, we have $i=j$ and $c(Z_i)=\Lambda$ and the stack of $N$ can only be empty. We want to illustrate the evaluation of a formula as follows where the stack is given by $ c(Z_i)\,c(Z_{i-1})\,\cdots\,c(Z_{18})$ ($i\geq 18$) in our representation. 

\vspace{0.2cm} 

When a new symbol is stored in the stack, this current symbol is displayed in \textcolor{blue}{blue}. The last symbol at the right end of the rest of the inputted formula currently read is displayed in \textcolor{red}{red}. The individuals $\land_0,\lor_0,\leftrightarrow_0, \land_1,\lor_1,\leftrightarrow_1,q_0$, and $1$ can be stored in $Z_{35}$ and used as symbols for denoting the internal state of the pushdown automaton $N$. At any time, the further simulation of the pushdown automaton depends on this state, the last symbol (at the left-hand end) of the stack and the last symbol (at the right-hand end) of the rest of the inputted formula. This means that three symbols are to be evaluated by the use of branching instructions. $c(Z_{35})=1$ means here that the input can be accepted. The change of the state from $q_0$ to $q_0$ means that a symbol $\textcolor{blue}{A}$ is simultaneously stored in the stack if the read symbol of $\phi$ was this symbol $\textcolor{red}{A}$ or if it was $\textcolor{red}{\neg}$ and the last symbol in the stack was a $B$ such that $\textcolor{blue}{A}=\textcolor{red}{\neg} B$ holds. Simultaneously, the red symbol $\textcolor{red}{A}$ and $\textcolor{red}{\neg}$, respectively, is removed such that we get a new rest of the formula. The change from $q_0$ to $\leftrightarrow_1$ means that the symbol at the right end of $\phi$ was $\textcolor{red}{\leftrightarrow}$ and the last stack symbol was $1$. Both symbols, $ \textcolor{red}{\leftrightarrow}$ and the last symbol 1 in the stack, are removed. The change from $\leftrightarrow_1$ to $q_0$ means that a symbol $\textcolor{blue}{A}$ is simultaneously stored in the stack if the last symbol in the stack was an $B$ such that $A= (B \leftrightarrow 1)$ holds, and so on. Each of the following symbols is stored in one of the $Z$-registers $Z_1,Z_2,\ldots$ for evaluating a symbol of a formula, a symbol stored in the stack, and an internal state of $N$. For deleting a symbol $c(Z_i)$ stored in the stack, we use an index register. If $c(I_k)=i$ holds, then the execution of the pseudo instruction $I_k:=I_k-1$ can have the same effect as deleting the individual $c(Z_i)$.

\vspace{0.2cm}

\begin{tabular}{lrc} 

$c(Z_1)\,\, c(Z_2)\,\,\quad\cdots \cdots \quad \,\,c(Z_{17})$&$ c(Z_{34})\,\,\,\,\,\cdots \,\,\,\,\,c(Z_{18})$&$\quad c(Z_{35})$\vspace{0.1cm}\\
$\underbrace{\land \lor \land 01 \neg\leftrightarrow 10 \lor \neg \land 01 \leftrightarrow 0\textcolor{red}{0} } _ {\phi {\rm \, after\, the \,input\,}}$ \hspace*{0.8cm} &$ \underbrace{\textcolor{blue}{\hspace{2.8cm}\#}} _ {{\rm the \,stack\,} \, {\rm at\, the\, beginning}} $&$\underbrace{ q_0}_{{\rm the \,state}}$\\
\end{tabular} 

\vspace{0.6cm}

The evaluation:

\begin{tabular}{lrc} 

$\overbrace{\land \lor \land 01 \neg\leftrightarrow 10 \lor \neg \land 01 \leftrightarrow 0\textcolor{red}{0} } ^ {\phi }$\hspace*{1.8cm}&$ \overbrace{\textcolor{blue}{\hspace{1.8cm}\#}} ^ {{\rm the \,stack\,} \,} $&$\overbrace{ q_0}^{{\rm the \,state}}$\\

$\land \lor \land 01 \neg\leftrightarrow 10 \lor \neg \land 01 \leftrightarrow \textcolor{red}{0} $ & $ \textcolor{blue}{0}\#$ & $q_0$\\

$\land \lor \land 01 \neg\leftrightarrow 10 \lor \neg \land 01 \textcolor{red}{ \, \leftrightarrow} $ & $ \textcolor{blue}{0}0\#$ & $q_0$\\

$\land \lor \land 01 \neg\leftrightarrow 10 \lor \neg \land 01 \textcolor{red}{ \, } $ & $ \textcolor{blue}{}0\#$ & $\leftrightarrow_0$\\

$\land \lor \land 01 \neg\leftrightarrow 10 \lor \neg \land 0 \textcolor{red}{ 1 } $ & $ \textcolor{blue}{1}\#$ & $q_0$\\

$\land \lor \land 01 \neg\leftrightarrow 10 \lor \neg \land \textcolor{red}{ 0 } $ & $ \textcolor{blue}{1}1\#$ & $q_0$\\

$\land \lor \land 01 \neg\leftrightarrow 10 \lor \neg \textcolor{red}{\,\land } $ & $ \textcolor{blue}{0}11\#$ & $q_0$\\

$\land \lor \land 01 \neg\leftrightarrow 10 \lor \neg \textcolor{red}{ } $ & $ \textcolor{blue}{}11\#$ & $\land_0$\\

$\land \lor \land 01 \neg\leftrightarrow 10 \lor \textcolor{red}{ \neg } $ & $ \textcolor{blue}{0}1\#$ & $q_0$\\

$\land \lor \land 01 \neg\leftrightarrow 10 \textcolor{red}{\, \lor } $ & $ \textcolor{blue}{1}1\#$ & $q_0$\\

$\land \lor \land 01 \neg\leftrightarrow 10 \textcolor{red}{\, } $ & $ \textcolor{blue}{}1\#$ & $\lor_1 $\\

$\land \lor \land 01 \neg\leftrightarrow 1 \textcolor{red}{0 } $ & $ \textcolor{blue}{1}\#$ & $q_0$\\

$\land \lor \land 01 \neg\leftrightarrow \textcolor{red}{1 } $ & $ \textcolor{blue}{0}1\#$ & $q_0$\\

$\land \lor \land 01 \neg \textcolor{red}{\,\leftrightarrow } $ & $ \textcolor{blue}{1}01\#$ & $q_0$\\

$\land \lor \land 01 \neg \textcolor{red}{ } $ & $ \textcolor{blue}{}01\#$ & $\leftrightarrow_1$\\

$\land \lor \land 01 \textcolor{red}{\neg } $ & $ \textcolor{blue}{0}1\#$ & $q_0$\\

$\land \lor \land 0 \textcolor{red}{1 } $ & $ \textcolor{blue}{1}1\#$ & $q_0$\\

$\land \lor \land \textcolor{red}{0 } $ & $ \textcolor{blue}{1}11\#$ & $q_0$\\

$\land \lor \textcolor{red}{\land } $ & $ \textcolor{blue}{0}111\#$ & $q_0$\\

$\land \lor \textcolor{red}{ } $ & $ \textcolor{blue}{}111\#$ & $\land _0$\\

$\land \textcolor{red}{ \lor } $ & $ \textcolor{blue}{0}11\#$ & $q _0$\\

$\land \textcolor{red}{ } $ & $ \textcolor{blue}{}11\#$ & $\lor _0$\\

$ \textcolor{red}{ \land } $ & $ \textcolor{blue}{1}1\#$ & $q _0$\\

$ \Lambda$ & $ \textcolor{blue}{}1\#$ & $\land _1$\\

$ \Lambda$ & $\textcolor{blue}{1}\#$ & $q_0$\\

$ \Lambda$ & $\#$ & $1$\\

$ \Lambda$ & $\Lambda $ & $1$\\\\
\end{tabular}

The machine halts only if the internal state of the automaton simulated is $1$. 
\end{example}

\noindent It is also possible to check the truth values of all quantifier-free first-order formulas by a single ${\cal A}$-machine if the codes of the machines described in Example \ref{Example2} are used because a universal machine (see \cite{GASS20}) can simulate the machines each of which allows to evaluate a single formula. The satisfiability of these formulas can be checked by non-deterministic (universal) machines (cf. Example \ref{PositivRaten}). Non-deterministic oracle machines can be used for evaluating first-order formulas with any number of quantifiers.

\vspace{0.1cm}

\section*{Summary}

 The {\em execution of an instruction} by an ${\cal A}$-machine ${\cal M}$ can be described purely mathematically by functions in ${\cal F}_{\cal M}$ (cf. also \cite{GASS20}). 

\begin{overview}
[Instructions and the change of configurations I]

\hfill

\nopagebreak 

\noindent \fbox{\parbox{11.8cm}{\small \em\begin{tabular}{l}

(1) ${\rm F}$-instructions { $\ell \!: \, Z_j:= f_i^{m_i}(Z_{j_1},\ldots, Z_{j_{m_i}}) $} \hfill $ _{(\Rightarrow\ell\in{\cal L}_{{\cal M},{\rm F}})}$\vspace{0.2cm}\\

\qquad $(\ell \,.\,\vec \nu\,.\,\bar u) \textcolor{red}{ \,\to_{\cal M}\,} (\ell +1 \,.\,\vec \nu\,.\, \underbrace {( u_1,\ldots,u_{j-1},f_i(u_{j_1},\ldots, u_{j_{m_i}}),u_{j+1}, \ldots)}_{ F_{\ell }(\bar u)}) $ 

\vspace{0.3cm}\\

(2) ${\rm F}_0 $-instructions { $\ell \!: \ Z_j:= c_i^0 $} \hfill $ _{(\Rightarrow\ell\in{\cal L}_{{\cal M},{\rm F}_0})}$\vspace{0.2cm}\\

\qquad $(\ell \,.\,\vec \nu\,.\,\bar u)\textcolor{red}{ \,\to_{\cal M}\,} (\ell +1 \,.\,\vec \nu\,.\, \underbrace {( u_1,\ldots,u_{j-1},c_{\alpha_i},u_{j+1}, \ldots) }_{ F_{\ell }(\bar u)})$

\vspace{0.3cm}\\

(3) ${\rm C}$-instructions { $\ell \!: \,Z_{I_j}:=Z_{I_k}$} \hfill $ _{(\Rightarrow\ell\in{\cal L}_{{\cal M},{\rm C}})}$\vspace{0.2cm}\\

\qquad $(\ell \,.\,\vec \nu\,.\,\bar u)\textcolor{red}{ \,\to_{\cal M}\,} (\ell +1 \,.\,\vec \nu\,.\,\underbrace {(u_1,\ldots,u_{\nu_j-1},u_{\nu_k},u_{\nu_j+1}, \ldots)}_{\ C_{\ell }(\vec \nu,\bar u)} )$

\vspace{0.3cm}\\

(4) ${\rm T}$-instructions 

{ \sf $\ell \!: \,$ if $r_i^{k_i}(Z_{j_1},\ldots, Z_{j_{k_i}})$ then goto $\ell _1$ else goto 
 $\ell _2$}\hspace{0.37cm}$_{(\Rightarrow\ell \in {\cal L}_{{\cal M},{\rm T}})}$\vspace{0.2cm}\\

\qquad $(\ell \,.\,\vec \nu\,.\,\bar u)\textcolor{red}{ \,\to_{\cal M}\,} (\underbrace{\ell_1}_{ T_{\ell }(\bar u)} .\,\vec \nu\,.\,\bar u)$ \qquad\qquad if $ (u_{j_1},\ldots, u_{j_{k_i}})\in r_i$ \\\

\qquad $(\ell \,.\,\vec \nu\,.\,\bar u) \textcolor{red}{\, \to_{\cal M}\,} (\,\overbrace{\ell_2}^{}.\,\vec \nu\,.\,\bar u)$ \qquad\qquad if $ (u_{j_1},\ldots, u_{j_{k_i}})\not \in r_i$

\vspace{0.3cm}\\

(5) ${\rm H}_{\rm T}$-instructions { \sf $\ell \!: \,$ if $I_j=I_k$ then goto $\ell _1$ else goto 
 $\ell _2$}\quad \hfill $_{(\Rightarrow \ell \in {\cal L}_{{\cal M},{\rm H}_{\rm T}})}$\vspace{0.2cm}\\

\qquad $(\ell \,.\,\vec \nu\,.\,\bar u) \textcolor{red}{ \,\to_{\cal M}\,} (\underbrace{\ell_1}_{ T_{\ell }(\vec \nu)} .\,\vec \nu\,.\,\bar u)$ \qquad\qquad if $ \nu_j=\nu_k$ \\

\qquad $(\ell \,.\,\vec \nu\,.\,\bar u)\textcolor{red}{ \,\to_{\cal M}\,} (\,\overbrace{\ell_2}^{}.\,\vec \nu\,.\,\bar u)$ \qquad\qquad if $\nu_j\not = \nu_k$

\vspace{0.3cm}\\

(6) ${{\rm H}_1}$-instructions { $\ell \!: \, I_j:=1$} \hfill $ _{(\Rightarrow\ell\in{\cal L}_{{\cal M},{\rm H}_1})}$\vspace{0.2cm}\\

\qquad $(\ell \,.\,\vec \nu\,.\,\bar u) \textcolor{red}{ \,\to_{\cal M}\,} (\ell +1\,.\, \underbrace {( \nu_1,\ldots,\nu_{j-1},1,\nu_{j+1}, \ldots)}_{ H_{\ell }(\vec \nu)}\,.\,\bar u )$ 
\vspace{0.3cm}\\

(7) ${{\rm H}_{+1}}$-instructions { $\ell \!: \, I_j:=I_j+1$} \hfill $ _{(\Rightarrow\ell\in{\cal L}_{{\cal M},{\rm H}_{+1}})}$\vspace{0.2cm}\\

\qquad $(\ell \,.\,\vec \nu\,.\,\bar u) \textcolor{red}{ \,\to_{\cal M}\,} (\ell +1\,.\, \underbrace {( \nu_1,\ldots,\nu_{j-1},\nu_j+1,\nu_{j+1}, \ldots)}_{ H_{\ell }(\vec \nu)}\,.\,\bar u )$ 
\vspace{0.3cm}\\

(8) ${\rm S}$-instruction { $\ell _{\cal P}\!: $ {\sf stop}} \hfill $ _{(\Rightarrow\ell_{\cal P}\in{\cal L}_{{\cal M},{\rm S}})}$\vspace{0.2cm}\\

\qquad $(\ell_{\cal P}\,.\,\vec \nu\,.\,\bar u) \textcolor{red}{ \,\to_{\cal M}\,} (\ell_{\cal P}\,.\,\vec \nu\,.\,\bar u) $

\end{tabular}}}
\end{overview}

\begin{overview}
[Instructions and the change of configurations II]

\hfill

\nopagebreak 

\noindent \fbox{\parbox{11.8cm}{\small \em\begin{tabular}{l}

(9) ${\rm O}$-instructions 
 \hspace{0.37cm}{\sf $\ell :\,$ if $(Z_1,\ldots,Z_{I_1})\in \!{ \cal O}$ then goto $\ell_1$ else goto $\ell_2$}\,\,\,$_{(\Rightarrow\ell \in {\cal L}_{{\cal M},{\rm O}})}$\vspace{0.2cm}\\

\qquad $(\ell \,.\,\vec \nu\,.\,\bar u)\textcolor{red}{ \,\to_{\cal M}\,} (\underbrace{\ell_1}_{ T_{\ell }(\nu_1,\bar u)} .\,\vec \nu\,.\,\bar u)$ \qquad\qquad if $ (u_1,\ldots, u_{\nu_1})\in Q$ \\

\qquad $(\ell \,.\,\vec \nu\,.\,\bar u) \textcolor{red}{\, \to_{\cal M}\,} (\,\,\,\overbrace{\,\,\ell_2\,\,}^{}\,\,\,\,.\,\vec \nu\,.\,\bar u)$ \qquad\qquad if $ (u_1,\ldots, u_{\nu_1})\not\in Q$

\vspace{0.4cm}\\

(10) ${\rm N}$-instructions \quad {$\ell:\,\, Z_j:=\nu[{\cal O}](Z_1,\ldots,Z_{I_1})$} \hfill $ _{(\Rightarrow\ell\in{\cal L}_{{\cal M},{\rm N}})}$\vspace{0.2cm}\\

\qquad $((\ell \,.\,\vec \nu\,.\,\bar u) , (\ell +1 \,.\,\vec \nu\,.\,\underbrace{( u_1,\ldots,u_{j-1}, \underbrace {\qquad y\qquad}_{\in \nu[Q] (u_1,\ldots, u_{\nu_1})},u_{j+1}, \ldots)}_{((\nu_1,\bar u), \ldots)\in F_\ell\quad}))\in\textcolor{red}{ \,\genfrac{}{}{0pt}{2}{\longrightarrow} {\longrightarrow}_{\cal M}\,} $ 

\vspace{0.4cm}\\

(11) ${\rm B}$-instructions 

{ \sf $\ell \!: \, $ goto $\ell _1$ else goto 
 $\ell _2$}\hfill $_{(\Rightarrow\ell \in {\cal L}_{{\cal M},{\rm B}})}$\vspace{0.2cm}\\

\qquad $((\ell \,.\,\vec \nu\,.\,\bar u), (\underbrace{\ell_0}_{\,\,\,(\ell,\ldots) \in G} .\,\vec \nu\,.\,\bar u))\in\textcolor{red}{ \,\genfrac{}{}{0pt}{2}{\longrightarrow} {\longrightarrow}_{\cal M}\,}$ 

\qquad \\ \end{tabular}}}
\end{overview}

\begin{overview}
[Non-deterministic components and instructions]

\hfill

\nopagebreak 

\noindent \fbox{\parbox{11.8cm}{

{\em Multi-valued input functions}

\vspace{0.1cm}

\quad without restriction

\qquad ${\rm In}_{\cal M}=\{(\vec x,((n, 1,\ldots, 1)\,.\,\vec x\,.\,\vec y\,.\,(x_n,x_n, \ldots) ))\mid \vec x\in U_{\cal A}^\infty\,\,\&\,\,\vec y \in U_{\cal A}^\infty\}$

\vspace{0.1cm}

\quad limited to two options for each of the guesses

\qquad ${\rm In}_{\cal M}=\{(\vec x,((n, 1,\ldots, 1)\,.\,\vec x\,.\,\vec y\,.\,(x_n,x_n, \ldots) ))\mid \vec x\in U_{\cal A}^\infty\,\,\&\,\,\vec y \in \{c_1,c_2\}^\infty\}$

\hfill for constants $c_1,c_2\in U_{\cal A}$

\vspace{0.1cm}

{\em Non-deterministic instructions}

\qquad $\ell:\,\, Z_j:=\nu[{\cal O}](Z_1,\ldots,Z_{I_1})$

\qquad {\sf $\ell :\,$ goto $\ell_1$ or goto $\ell_2$}

\vspace{0.1cm}
}}
\end{overview}

\begin{overview}
[Several result functions]

\hfill

\nopagebreak 

\noindent \fbox{\parbox{11.8cm}{

{\em Deterministically computed result functions \textcolor{blue}{${\rm Res}_{\cal M}: \,\subseteq U_{\cal A}^\infty \to U_{\cal A}^\infty $}}

\vspace{0.1cm}

\qquad ${\rm {\rm Res}_{\cal M}}(\vec x)={\rm Output}_{\cal M}( (\to_{\cal M})_{{\rm Stop}_{\cal M}} ( {\rm Input}_{\cal M}({\vec x})))$ \hfill {\rm (\ref{Res1})}

\vspace{0.1cm}

{\em Non-deterministically computed result functions \textcolor{blue}{${\rm Res}_{\cal M}: U_{\cal A}^\infty\to {\mathfrak P}(U_{\cal A}^\infty)$}}

\vspace{0.1cm}

\quad resulting from multi-valued input functions

 \qquad ${\rm Res_{\cal M}}(\vec x)=\{ {\rm Ouput}_{\cal M}((\to_{\cal M})_{{\rm Stop}_{\cal M}} (con))\mid (\vec x, con)\in {\rm Input}_{\cal M}(\vec x)\}$ \hfill {\rm (\ref{Res2})}

\vspace{0.1cm}

\quad resulting from non-deterministic instructions

\qquad ${\rm {\rm Res}_{\cal M}}({\vec x})=\{{\rm Output}_{\cal M}(con)\mid ({\rm Input}_{\cal M}(\vec x),con)\in (\genfrac{}{}{0pt}{2}{\longrightarrow} {\longrightarrow}_{\cal M})_{{\rm Stop}_{\cal M}}\}$ \hfill {\rm (\ref{Res3})}

\vspace{0.1cm}
}}
\end{overview}

\section*{Outlook: Simulations and complete problems}\label{Outlook}
\addcontentsline{toc}{chapter}{\bf Outlook}
\markboth{OUTLOOK}{OUTLOOK}

\begin{overview}[Hierarchies and complete problems]\label{KnownOrExpected}\hfill

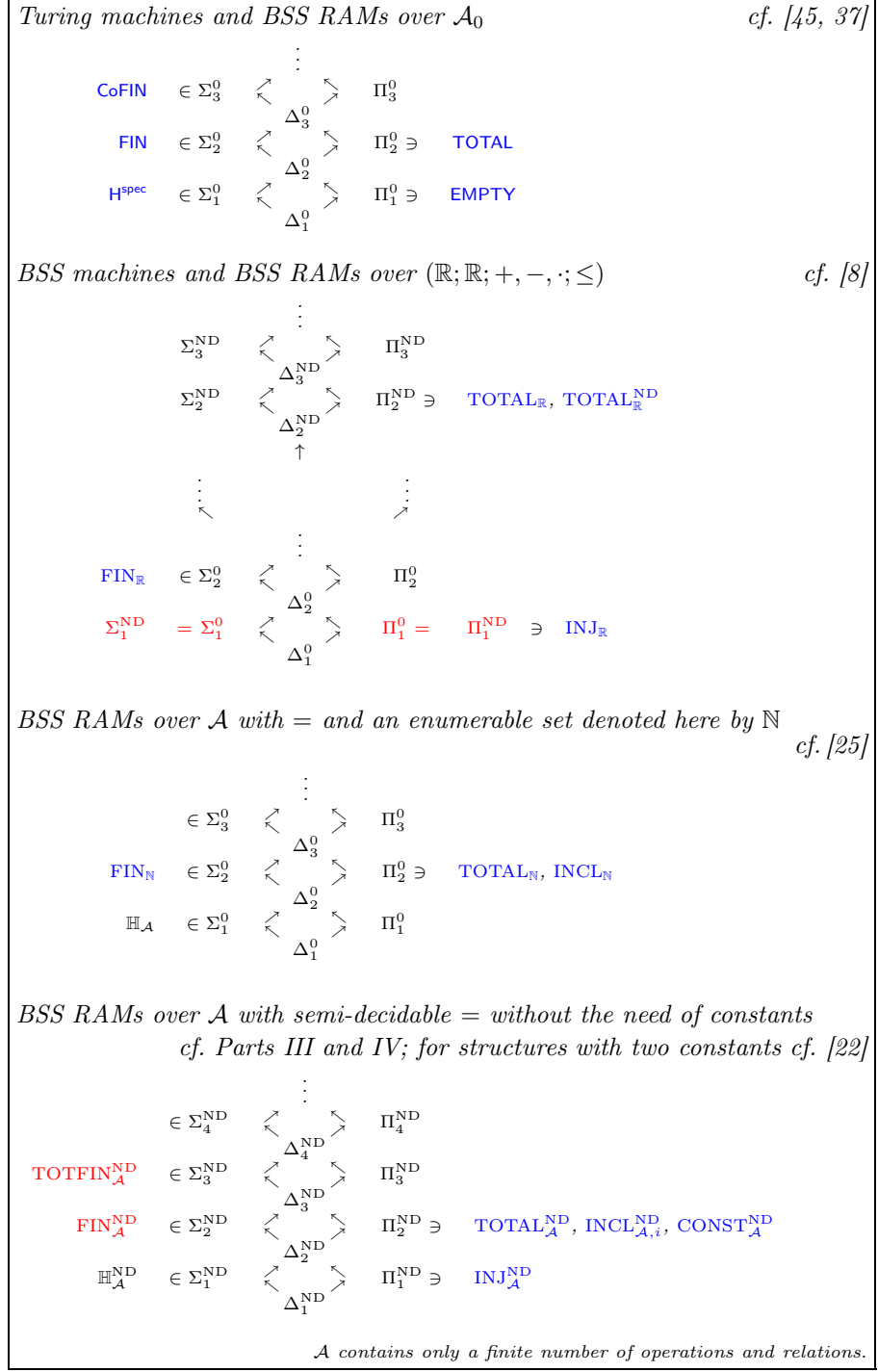
\nopagebreak

\noindent
\fbox{\parbox{11.75cm}
{\scriptsize
{\normalsize Turing machines and BSS RAMs over ${\cal A}_0$ \hfill cf.\,\,\cite{SOARE, KOZEN}}

\hspace{0.9cm}\begin{tabular}{rccccclcc}
&&$\vdots$&\\
\textcolor{blue}{$\mbox{\sf CoFIN}$}&$\in\Sigma_3^0 $&${\nearrow\atop\nwarrow}\qquad{\nwarrow\atop\nearrow}$&$\Pi_3^0\quad $&&\quad \\
&& $\Delta_3^0$\\
\textcolor{blue}{$\mbox{\sf FIN}$}&$\in\Sigma_2^0 $&${\nearrow\atop\nwarrow}\qquad{\nwarrow\atop\nearrow}$&$\Pi_2^0 \ni$&\textcolor{blue}{$\mbox{\sf TOTAL}$}&\\
&& $\Delta_2^0 $\\
\textcolor{blue}{$\mbox{\sf H}^{\sf spec}$}&$\in\Sigma_1^0 $&${\nearrow\atop\nwarrow}\qquad{\nwarrow\atop\nearrow}$&$\Pi_1^0\ni$&\textcolor{blue}{$\mbox{{\sf EMPTY}}$}\\
&& $\Delta_1 ^0$\\\\
\end{tabular}

\noindent {\normalsize BSS machines and BSS RAMs over $(\mathbb{R};\mathbb{R};+,-,\cdot;\leq)$ \hfill cf.\,\,\cite{CUCKER92}}

\hspace{0.95cm}\begin{tabular}{rccclcc}&&$\vdots$&\\
&$\Sigma_3^{\rm ND}$&${\nearrow\atop\nwarrow}\qquad{\nwarrow\atop\nearrow}$&$\Pi_3^{\rm ND}$\\
&& $\Delta_3^{\rm ND}$\\
&$\Sigma_2^{\rm ND}$&${\nearrow\atop\nwarrow}\qquad{\nwarrow\atop\nearrow}$&$\Pi_2^{\rm ND}\ni$&\textcolor{blue}{$\mbox{{\rm TOTAL}}_\mathbb{R}$}, $\textcolor{blue}{\mbox{\rm TOTAL}_\mathbb{R}^{\rm ND}}$\\
&& $\Delta_2^{\rm ND}$\\
&&$\uparrow$&\\
&${\vdots\,\,\,\,\atop\nwarrow}\!\!\!$& &$\!\!\!{\vdots\!\!\!\!\atop\nearrow}$&\\
&&$\vdots$&\\
\textcolor{blue}{$\mbox{\rm FIN}_\mathbb{R}$}&$\in\Sigma_2 ^0$&${\nearrow\atop\nwarrow}\qquad{\nwarrow\atop\nearrow}$&$\Pi_2^0 $\\
&& $\Delta_2^0 $\\
\textcolor{red}{$\Sigma_1^{\rm ND}$}&\textcolor{red}{$=$} \textcolor{red}{$\Sigma_1^0 $}&${\nearrow\atop\nwarrow}\qquad{\nwarrow\atop\nearrow}$&\textcolor{red}{$\Pi_1^0 $} \textcolor{red}{$=$}&\textcolor{red}{$ \Pi_1^{\rm ND}$}$\quad \ni$\quad \textcolor{blue}{$\mbox{{\rm INJ}}_\mathbb{R}$}\\\
&& $\Delta_1^0 $&&&&\\\\
\end{tabular}

\vspace{0.2cm}

\noindent {\normalsize BSS RAMs over ${\cal A}$ with $=$ and an enumerable set denoted here by $\bbbn$}

\hfill {\normalsize cf.\,\cite{GASS17}}

\hspace{1.08cm}\begin{tabular}{rccclcc}
&&$\vdots$&\\
&$\in\Sigma_3^0 $&${\nearrow\atop\nwarrow}\qquad{\nwarrow\atop\nearrow}$&$\Pi_3^0\quad $&&\quad \\
&& $\Delta_3^0$\\
\textcolor{blue}{$\mbox{\rm FIN}_\bbbn$}&$\in\Sigma_2^0 $&${\nearrow\atop\nwarrow}\qquad{\nwarrow\atop\nearrow}$&$\Pi_2^0 \ni$&\textcolor{blue}{$\mbox{{\rm TOTAL}}_\bbbn$}, \textcolor{blue}{$\mbox{{\rm INCL}}_\bbbn$}&&\\
&& $\Delta_2^0 $\\
$\bbbh_{\cal A}$&$\in\Sigma_1^0 $&${\nearrow\atop\nwarrow}\qquad{\nwarrow\atop\nearrow}$&$\Pi_1^0\quad $&&\quad \\
&& $\Delta_1 ^0$
\\\\
\end{tabular} 

\vspace{0.2cm}
 
\noindent {\normalsize BSS RAMs over ${\cal A}$ with semi-decidable $=$ without the need of constants} 

\vspace{0.1cm}

\hfill {\normalsize cf.\,\,Parts III and IV; for structures with two constants cf.\,\,\cite{GASS15}}

\begin{tabular}{rccclcc}
&&$\vdots$&\\
&$\in\Sigma_4^{\rm ND}$&${\nearrow\atop\nwarrow}\qquad{\nwarrow\atop\nearrow}$&$\Pi_4^{\rm ND}\quad $&&\quad \\
&& $\Delta_4^{\rm ND}$\\
\textcolor{red}{$\mbox{\rm TOTFIN}_{\cal A}^{\rm ND}$}&$\in\Sigma_3^{\rm ND}$&${\nearrow\atop\nwarrow}\qquad{\nwarrow\atop\nearrow}$&$\Pi_3^{\rm ND}\quad$\\
&& $\Delta_3^{\rm ND}$\\
\textcolor{red}{$\mbox{\rm FIN}_{\cal A}^{\rm ND}$}&$\in\Sigma_2^{\rm ND}$&${\nearrow\atop\nwarrow}\qquad{\nwarrow\atop\nearrow}$&$\Pi_2^{\rm ND}\ni$&\textcolor{blue}{$\mbox{{\rm TOTAL}}_{\cal A}^{\rm ND}$}, \textcolor{blue}{$\mbox{{\rm INCL}}_{{\cal A},i}^{\rm ND}$}, \textcolor{blue}{$\mbox{{\rm CONST}}_{\cal A}^{\rm ND}$}\\
&& $\Delta_2^{\rm ND}$\\
$\bbbh_{\cal A}^{\rm ND}$&$\in\Sigma_1^{\rm ND}$&${\nearrow\atop\nwarrow}\qquad{\nwarrow\atop\nearrow}$&$\Pi_1^{\rm ND}\ni$&\textcolor{blue}{$\mbox{{\rm INJ}}_{\cal A}^{\rm ND}$} &&\\
&& $\Delta_1^{\rm ND}$\\&&&&\\
\end{tabular}

\hfill ${\cal A}$ contains only a finite number of operations and relations.
}}
\end{overview}

\vspace{0.2cm}

\noindent
In Part II, we compare non-deterministic BSS RAMs for a given structure ${\cal A}$ and give a characterization of these machines by using Moschovakis' operator. In Part III, we will consider non-deterministic universal (oracle) machines and more general halting problems, such as the problems $^{\rm general}\mathbb{H}_{\cal A}^{\big[{\rm ND}\big]}$, for deterministic and non-deterministic machines and, in particular, for oracle machines that we will need for a characterization of hierarchies defined analogously to the arithmetical hierarchy. We think that Part II and Part III can help to better understand the similarities between the hierarchies introduced by Moschovakis in \cite{MOSCHO} and other hierarchies of undecidable decision problems. 

Based on studies in \cite{SOARE, KOZEN, MEERZIEGLER06, GASS08A,GASS13} and in \cite{CUCKER92} and in \cite{GASS17} we want to continue the classification of decision problems and consider entire hierarchies of classes of decision problems and complete problems for these classes.

In \cite{GASS17}, we considered a first hierarchy defined by deterministic machines over several structures ${\cal A}$ as follows.
 \begin{displaymath} \begin{array}{lcl}
{\cal A}\mbox{-}{\textcolor{blue}{\Sigma_0^0 }}&= &{\rm DEC}_{\cal A}
 \\{\cal A}\mbox{-}{\textcolor{blue}{\Pi_n ^0}}&=&\{U_ {\cal A}^\infty \setminus P\mid P \in \Sigma_n ^0\}\\
{\cal A}\mbox{-}{\textcolor{blue}{\Delta_n^0 }}&=& \Sigma_n^0 \cap \Pi_n ^0\vspace{0.1cm}\\
{\cal A}\mbox{-}{\textcolor{blue}{\Sigma_{n+1 } ^0}}&= &\{P\subseteq U_{\cal A}^\infty\mid (\exists Q\in\Pi_n ^0)(\textcolor{red}{ P \in {\rm SDEC}_{\cal A}^Q})\}
\end{array}\end{displaymath}
Analogously to the arithmetical hierarchy, a second hierarchy can be defined syntactically by formulas (of an infinitary logic):
\begin{displaymath} \begin{array}{lcl}
{\cal A}\mbox{-}{\textcolor{blue}{\Sigma_0^{\rm ND}}}&= &{\rm DEC}_{\cal A}
 \vspace{0.1cm} \\{\cal A}\mbox{-}{\textcolor{blue}{\Pi_n^{\rm ND}}}&=&\{ U_{\cal A}^\infty \setminus P\mid P \in \Sigma_n^{\rm ND}\} \vspace{0.1cm}\\
{\cal A}\mbox{-}{\textcolor{blue}{\Delta_n^{\rm ND}}}&=& \Sigma_n^{\rm ND} \cap \Pi_n^{\rm ND} \vspace{0.1cm}\\
{\cal A}\mbox{-}{\textcolor{blue}{\Sigma_{n+1 }^{\rm ND}}}&= &\{P\subseteq U_{\cal A}^\infty\mid (\exists Q\in\Pi_n^{\rm ND}) \textcolor{red}{\forall\vec x(\vec x\in P\Leftrightarrow (\exists \vec y \in U_{\cal A}^\infty)(\langle \vec x ,\vec y\rangle\in Q))} \}
\end{array}\end{displaymath}
where $\langle \vec x ,\vec y\rangle=_{\rm df}(x_1,\ldots,x_n,y_1,\ldots,y_m,x_n,\ldots,x_n)\in U_{\cal A}^{ca(n,m)}$. 
The resulting hierarchies are presented in Overview \ref{KnownOrExpected} where the prefix ${\cal A}\mbox{-}$ is omitted and the arrows $\to$ stand for {\it is included in}. The both first hierarchies for Turing machines and for BSS machines are discussed, for instance, in \cite{SOARE} and \cite{KOZEN} and studied in \cite{CUCKER92}, respectively. The third hierarchy was considered in \cite{GASS17}. The fourth hierarchy and some complete problems (written in blue for some of the classes of the latter hierarchy in Overview \ref{KnownOrExpected}) were presented and discussed for structures ${\cal A}$ containing at least two constants at the CCC 2015 (cf. \cite{GASS15}). Note, that the latter hierarchy can also be characterized by nondeterministic machines over ${\cal A}$ because we have
\begin{displaymath} \begin{array}{lcl}
{\cal A}\mbox{-}{\textcolor{blue}{\Sigma_{n+1 }^{\rm ND}}}&= &\{P\subseteq U_{\cal A}^\infty\mid (\exists Q\in\Pi_n )(\textcolor{red}{ P \in (\mbox{\rm SDEC}_{\cal A}^{\rm ND})^Q})\}.
\end{array}\end{displaymath}
In Part IV, we will define several complete problems for some of the classes of the latter hierarchy without the need to have constants and provide the necessary proofs.
While, in most cases, we considered structures with two constants that allow to use codes in ASCII format and the like, we will in the following use G\"odel numbers as codes stored in index registers (cf.\,\,Consequences \ref{SideEffect}). 

\vspace{0.3cm}

We would like to point out that we try to characterize the classical term {\sf algorithm} in a broad and comprehensive sense. To describe the concept of abstract computability, we use some terms from computer science, even if this does not mean that the considered programs can be executed by computers for the given first-order structures. This is only possible for certain structures.

\section*{Acknowledgment}
\addcontentsline{toc}{chapter}{\bf Acknowledgment}
\markboth{ACKNOWLEDGMENT}{ACKNOWLEDGMENT}

I would like to thank the participants of my lectures Theory of Abstract Computation and, in particular, Patrick Steinbrunner and Sebastian Bierba\ss{} for useful questions and the discussions. My thanks go also to all participants of meetings in Greifswald and Kloster. In particular, I would like to thank my co-authors Arno Pauly and Florian Steinberg for the discussions on operators related to several models of computation and Vasco Brattka, Philipp Schlicht, and Rupert H\"olzl for many interesting discussions.

{\small

\noindent For this article, we also used translators such as those of Google and DeepL. }

\vspace{4cm}

\vfill

\hfill {\scriptsize\sf C $\cdot $ H $\cdot$ T}
\end{document}